\documentclass[final,leqno]{article}

\usepackage{latexsym}
\usepackage{amsmath}%
\usepackage{amssymb}%

\newtheorem{remark}{Remark}
\newcommand\beq{\begin{equation}}
\newcommand\eeq{\end{equation}}

\def\Dbu{\boldsymbol{D}(\bu)}

\def\bn{\boldsymbol n}
\def\bN{\boldsymbol N}

\newcommand{\grad}{{\mathbf{\nabla}}}

\newcommand{\D}{\mathcal{D}}
\newcommand{\Wi}{{\text{Wi}}}
\newcommand{\Bi}{{\text{Bi}}}
\renewcommand{\Re}{{\text{Re}}}
\newcommand{\De}{{\text{De}}}
\newcommand{\e}{\varepsilon}
\newcommand{\R}{\mathbb{R}}
\newcommand{\f}{\boldsymbol{f}}
\newcommand{\I}{\boldsymbol{I}}
\newcommand{\xx}{\boldsymbol{x}}
\newcommand{\bolda}{\boldsymbol{a}}
\newcommand{\be}{\boldsymbol{e}}

\newcommand{\bzero}{{\boldsymbol{0}}}
\renewcommand{\div}{\operatorname{div}}
\newcommand{\curl}{\operatorname{curl}}
\newcommand{\tr}{\operatorname{tr}}
\newcommand{\Sb}{\boldsymbol{S}}
\newcommand{\Tb}{\boldsymbol{T}}

\newcommand{\gbu}{{\grad\bu}}
\newcommand{\bu}{{\boldsymbol{u}}}
\newcommand{\pds}[1]{\partial_{#1}}
\newcommand{\pdds}[2]{\partial^2_{#1#2}}
\newcommand{\str}{{\boldsymbol{\tau}}}
\newcommand{\strs}{{\boldsymbol{\sigma}}}
\newcommand{\bv}{\boldsymbol{v}}

\title{
Unified formal reduction for 
fluid models 
of free-surface shallow gravity-flows 
}

\author{
Fran\c{c}ois Bouchut
\thanks{
Universit\'e Paris-Est, Laboratoire d'Analyse et de Math\'ematiques Appliqu\'ees (UMR 8050), CNRS, UPEMLV, UPEC, F-77454, Marne-la-Vall\'ee, France
}
\and
S\'ebastien Boyaval
\thanks{
Universit\'e Paris-Est, Laboratoire d'hydraulique Saint-Venant, Ecole Nationale des Ponts et Chauss\'ees -- EDF R\&D -- CETMEF, 6 quai Watier, 78401 Chatou Cedex, France ; MICMAC team INRIA Rocquencourt ; 
sebastien.boyaval@enpc.fr (corresponding author)
}
}

\begin{document}

\maketitle

\begin{abstract}
We propose a unified 
approach to
the formal long-wave reduction of several fluid models for thin-layer incompressible homogeneous flows driven by a constant external force like gravity.
The procedure is based on a mathematical coherence property that univoquely defines one reduced model given one rheology and one thin-layer regime. For the first time, as far as we know, various known reduced models can thus be investigated within a single mathematical framework, for various rheologies (viscous and viscoelastic) and various limit regimes (fast inertial flows and slow viscous flows). Furthermore, our systematic procedure also produces 
new reduced models for viscoelastic non-Newtonian fluids and improves on our previous work [Bouchut \& Boyaval, M3AS (23) 8, 2013].
\end{abstract}

%

\pagestyle{myheadings}
\thispagestyle{plain}
\markboth{Bouchut and Boyaval}{Unified derivation of reduced models for shallow flows}

\section{Introduction}

Formal a priori simplifications of a model is a game physicists and mathematicians have been playing for years, in particular for fluid equations.
Historically, reduced models indeed proved useful because they were more amenable to analytical (exact) solution than full models, for instance the Saint-Venant equations, and thus helped understand a few simple phenomena (e.g. dam breaks) in a time where computer simulations did not exist.
Nowadays, reduced models are still sometimes preferred to full models, for instance to numerically simulate cheaply complex fluid flows, and next discriminate against various possible rheologies by comparison with experiments. In this work, we have more precisely in mind paving the way for a better modelling of the rheology in e.g. mud flows and landslides, which are still much investigated, by experimentalists in particular~\cite{ancey-2007}. Indeed, their thin-layer geometry apparently suits well with simplifications, see e.g.~\cite{ancey-cochard-2009}.

Numerous reduced models have already been derived for gravity-driven free-surface shallow flows with various rheologies in various regimes, however they are still difficult to connect one another. A unifying viewpoint of various surface wave models for water (Newtonian) flows has been constructed recently~\cite{bonneton-lannes-2009}, but it holds only for purely-irrotational water flow models. A generic mathematically-inclined derivation procedure for slow flows closed to a stationary solution has also been used recently for various rheologies, see e.g.~\cite{ruyerquil-manneville-2000,fernandez-nieto-noble-vila-2010}, 
but it holds only for viscous (laminar) flow regimes. Our primary goal here is to establish another mathematically-inclined framework, that is common to various {\it thin-layer} flows (slow/fast) and several fluids (Newtonian/non-Newtonian). To this aim, we introduce a procedure that yields various long-wave reduced models when using various rheologies with Navier-Stokes equations depending on various thin-layer flow regimes.

Our procedure, inspired by~\cite{gerbeau-perthame-2001,marche-2007} where the viscous shallow water equations are derived (see Section~\ref{sec:thin}), is based on a 
very natural mathematical coherence property.
Although it 
cannot certify {\it rigorously} when and how a solution to the reduced model is a good approximation of a solution to the full model, our formal procedure is univoque at least: it delivers a single reduced model given one rheology and one thin-layer long-wave flow regime,
whose solution should approximately solve the original full model provided the assumptions used for the derivation hold.

We obtain a synthetic viewpoint of various possible simplifications of the Navier-Stokes equations modelling a fluid under gravity in a free-surface thin-layer geometry, when the rheology varies (that is, the modelling of the internal stresses), as well as when the scaling assumptions for the flow regime vary (i.e. in the momentum equation, the hydrostatic forces are mainly balanced by the purely kinematic hydrodynamic forces in fast inertial flows, as opposed to the internal stresses in slow viscous flows). 

We also obtain new reduced models for fluids with complex rheologies.

\begin{itemize}
 \item For viscous Newtonian fluids (modelled by the standard Navier-Stokes equations), we obtain either viscous shallow water  equations in inertial regimes (as e.g. in~\cite{gerbeau-perthame-2001,marche-2007}) or lubrication equations in viscous regimes (as e.g. in~\cite{RevModPhys.69.931,ruyerquil-manneville-2000,matar-kraster-2009}), in Section~\ref{sec:newtonian}. 
 \item For viscous non-Newtonian fluids (nonlinear power-law models), we obtain either a nonlinear version of the shallow water equations
 in inertial regimes that is apparently new, or nonlinear lubrication equations in viscous regimes (see~\cite{fernandez-nieto-noble-vila-2010} and references therein), in Section~\ref{sec:power}.
 \item For viscoelastic non-Newtonian fluids, we obtain either shallow water equations with additional stress terms which extends the recent work~\cite{bouchut-boyaval-2013} in inertial regimes, or new lubrication equations in viscous regimes (different than those in~\cite{chupin-2009a,chupin-2009b}),  in Section~\ref{sec:viscoelastic}.
\end{itemize}

A few remarks are also in order.
\begin{itemize}
 \item The case of perfect fluids (no internal stresses) is singular, and it seems one cannot naturally derive a closed reduced model without coming back to an assumption about the dissipation terms initially neglected: we treat it as the inviscid limit of the viscous Newtonian case.
 \item The case of some viscoplastic Non-Newtonian fluids (i.e. some Bingham fluids) occurs as a singular limit of the nonlinear power-law models. Now, this case is interesting from the modelling viewpoint (some plastic non-Newtonian fluids are believed to possess a yield-stress that seems to suit well for modelling transitions of fluid-solid type like e.g. in avalanches), but already difficult from the mathematical viewpoint (the model is undetermined below the yield-stress) as well as from the physical viewpoint (the existence of a yield-stress to account for a transition of the fluid-solid type is still much debated).
 That is why it is in fact the single fluid model with plasticity that we investigate here. Our framework could nevertheless serve as a basis for future thin-layer investigations of models taking into account transitions of fluid-solid type, the modelling of viscoplasticity still being in its infancy.
 \item Concerning viscoelastic fluids, we improve here the model derived in~\cite{bouchut-boyaval-2013} from simple constitutive equations (only linear in the tensor state variable ``conformation'' that accounts for viscoelasticity, though already physically-consistent from the frame-invariance viewpoint) in the sense that here, we take into account friction at the bottom, surface tension, two-dimensional effects and a purely Newtonian additional viscosity. 
\end{itemize}
For a recent physically-inclined review of thin-film flows, we recommend~\cite{matar-kraster-2009}, and~\cite{RevModPhys.69.931} for an older one with a focus on stability. Let us now mathematically set the problem.

\section{Mathematical setting of the problem}

We endow the space $\R^3$ with a Galilean reference frame using cartesian coordinates $(\be_x,\be_y,\be_z)$.
We denote by $a_x$ (respectively $a_y$, $a_z$) the component in direction $\be_x$ (resp. $\be_y$, $\be_z$) of a vector (that is a rank-1 tensor) $\bolda$, by $a_{xx}, a_{xz},\ldots$ the components of higher-rank tensors, by $\bolda_H$ the vector of ``horizontal'' components $(a_x,a_y)$, by $(\bolda_H)^\bot=(-a_y,a_x)$ an orthogonal vector, by $\grad_H a$ the horizontal gradient $(\partial_xa,\partial_ya)$ of a smooth function $a:(x,y)\to a(x,y)$, and by $D_ta$ the material time-derivative $\pds{t}a + (\bu\cdot\grad)a$. 
We use the Frobenius norm $|\bolda|=\tr(\bolda^T\bolda)^{1/2}$ for a 2-tensor.


Gravity flows of incompressible homogeneous fluids are governed by Navier-Stokes equations
\begin{equation}
\label{eq:navier-stokes1}
\left\lbrace
\begin{aligned}
D_t\bu & = \div(\Sb) + \f & \text{ in } \D(t) \,,
\\
\div\bu & = 0 & \text{ in } \D(t) \,,
\end{aligned}
\right.
\end{equation}
on a regional scale where gravity is the single external force $\f$, with the velocity field $\bu$ as unknow variable plus
Cauchy stress tensor $\Sb = -p\I+\Tb$ ($\Tb$ is the deviatoric part of $\Sb$ when $\tr(\Tb)=0$).

We consider cases where the fluid is contained for all times $t\ge0$ within a cylindrical domain
\beq
\label{cylindrical_domain}
\D(t)= \{\xx=(x,y,z)\,,\ (x,y)\in\Omega_0\,, 0<z-b(x,y)<h(t,x,y) \} \,,
\eeq
with a {\it free surface} $z=b(x,y)+h(t,x,y)$ (a simplified modelling for a liquid-gas interface).
The free surface and the bottom topography $z=b(x,y)$ are thus {\sl unfolded} (i.e. single-valued) two-dimensional parametrized manifolds.
Furthermore, we are more specifically interested in the case of shallow flows where the two manifolds are assumed close to one-another in the ``vertical'' direction $\be_z$ -- whatever the ``horizontal'' position $(x,y)\in\Omega_0$ -- in comparison with a characateristic horizontal length $L$, and where they a priori never touch (though one usually next extends the application of the model to cases with vacuum).
We write this assumption
$$ h\sim\e $$
using an adimensional small parameter $\e\ge0$. It means that $h/(\e L)$ is bounded above and below independently of $(x,y)\in\Omega_0$, $t\ge 0$ as $\e\to0$,
as opposed to $a=O(\e)$ for a variable $a$, which simply means that the adimensional quantity $a/(\e A)$ (where $A$ is the natural characteristic size of $a$ as a function of $L,T$, see below) is bounded above, and may in fact decay faster than $\e$ to zero as $\e\to0$.
Note that we shall also use componentwise notation, e.g. $a_1,a_2=O(\e^{\alpha_1},\e^{\alpha_2})$ for $a_1=O(\e^{\alpha_1})$ and $a_2=O(\e^{\alpha_2})$.

The goal of this work is to derive a closed system of approximate equations for the flow 
that hopefully define a simpler and useful (that is, physically meaningful) mathematical model than~\eqref{eq:navier-stokes1} in the vicinity of the limit $\e\to0$. 
Note that we limit to cases where $\grad_H h = O(\e)$ holds, so that the reduced model captures only {\it long-wave} oscillations of the free surface.

Of course, at this stage, the system of equations~\eqref{eq:navier-stokes1} is not closed.
One still need to specify the rheology of the fluid (that is, invoke other equations linking $\Sb$ with $\bu$) as well as boundary conditions. 
We recall that it is exactly the goal of this work to derive approximations of~\eqref{eq:navier-stokes1} for {\it various} rheologies and flow regimes using the same procedure, thereby defining a common framework to compare rheological models with experimental measures through simple reduced models in the case of free-surface thin-layer flows. We have in mind rheological models for:
\begin{itemize}
 \item viscous Newtonian fluids like water, such that the deviatoric stress tensor is a linear function of the rate-of-strain tensor $\Dbu=\frac12(\gbu+\gbu^T)$, hence $\Tb=2\eta_s\Dbu$, with $\eta_s$ a {\it constant} kinematic viscosity,
 \item viscous non-Newtonian fluids, what most complex fluids are in a small range of shear rates at least, such that the mechanical behaviour is still decsribed with a purely viscous deviatoric stress tensor, but using a nonlinear power-law $\Tb=2\eta_s|\Dbu|^{n-1}\Dbu$ (termed pseudoplastic or shear-thinning if $0< n< 1$, dilatant or shear-thickening if $n>1$) for viscosity,
 \item viscoelastic non-Newtonian fluids like polymer solutions, such that $\Tb=2\eta_s\Dbu+\str$ invokes a non-Newtonian extra-stress $\str$ that is not necessarily deviatoric and defined through supplementary (integro-)differential equations.
\end{itemize}
Note that there is a huge amount of non-Newtonian models in the literature~\cite{barnes-hutton-walters-1989}.
\begin{itemize}
 \item Interestingly, the nonlinear power-law models for viscous non-Newtonian fluids coincide with the standard Navier-Stokes equations for viscous Newtonian fluids when $n=1$ and with a Bingham model for viscoplastic fluids when $n=0$. Although stress is undetermined in Bingham model when $|\Dbu|=0\Leftrightarrow|\Tb|<2\eta_s$, Bingham model can be understood as the limit of a regularized model~\cite{duvaut-lions-1972} and remains the least disputed basic model for the still much debated viscoplastic non-Newtonian fluids.
 \item Viscoelastic fluid models have been used successfully for the accurate description of polymer solutions for instance, see e.g.~\cite{bird-curtiss-armstrong-hassager-1987a,bird-curtiss-armstrong-hassager-1987b}. We repeat that we shall be content here with simple prototypical models among the numerous possibilities (see Section~\ref{sec:viscoelastic}).
\end{itemize}

We believe that the various prototypical rheologies mentionned above are representative enough to define a first common mathematical framework for various shallow flows with various rheologies.
We now complement them with the following common boundary conditions (BCs).

Let us denote by $\bn:(x,y)\in\Omega_0\to\bn(x,y)$ the unit vector of the direction normal to the bottom
\begin{equation}
\bn = \begin{pmatrix}
       - \grad_H b \\ 1
      \end{pmatrix} / \sqrt{1+|\grad_H b|^2}
\label{eq:normal_bottom}
\end{equation}
(inward the fluid) and by $(N_t,\bN)$ the time-space normal at the free surface (outward the fluid)
\begin{equation}
N_t = -\pds{t}{(b+h)}/\sqrt{1+|\grad_H(b+h)|^2}
 \qquad
\bN = \begin{pmatrix}
       - \grad_H (b+h) \\ 1
      \end{pmatrix} / \sqrt{1+|\grad_H (b+h)|^2} \,.
\label{eq:normal_top}
\end{equation}
An orthonormal frame is defined locally everywhere on the bottom using as basis in tangent planes
\begin{equation}
\label{eq:tangent}
\boldsymbol{t}_1 = 
\begin{pmatrix}
(\grad_H b)^\bot \\ 0 
\end{pmatrix}/|\grad_H b|
 \qquad
\boldsymbol{t}_2 = 
\begin{pmatrix}
- \grad_H b \\ -|\grad_H b |^2
\end{pmatrix}/(|\grad_H b |\sqrt{1+|\grad_H b |^2})
\end{equation}
when $|\grad_H b|\neq0$, otherwise $\boldsymbol{t}_1=(0,-1,0)^T$, $\boldsymbol{t}_2=(-1,0,0)^T$.
We require, at the bottom of the fluid, no penetration in the normal direction
\begin{equation}
\bu\cdot\bn=0  \,,
\quad\mbox{for }z=b(x,y),\ (x,y)\in\Omega_0\,,
\label{eq:nopenetration}
\end{equation}
and a Navier friction dynamic condition with coefficient $k$ in the tangent plane
\begin{equation}
\Sb\bn\wedge\bn = k \bu\wedge\bn \,,\quad\mbox{for }z=b(x,y),\ (x,y)\in\Omega_0\,;
\label{eq:friction}
\end{equation}
at the free surface, the usual kinematic condition $N_t+\bN\cdot\bu=0$, i.e.
\begin{equation}
-\partial_t(b+h) - u_x\partial_x(b+h) - u_y\partial_y(b+h) + u_z = 0\,, \quad\mbox{for }z=b(x,y)+ h(t,x,y),\ (x,y)\in\Omega_0\,,
\label{eq:kinematic}
\end{equation}
and surface tension with coefficient $\gamma$ as dynamic condition
\begin{equation}
\Sb\bN=\gamma\kappa\bN,
\quad\mbox{for }z=b(x,y)+ h(t,x,y),\ (x,y)\in\Omega_0,
\label{eq:tension}
\end{equation}
where $\kappa(t,x,y)=-\div\bN(t,x,y)$ is the (local) mean curvature at $z=b(x,y)+ h(t,x,y),\ (x,y)\in\Omega_0$; and finally, at the lateral boundary $\{\xx=(x,y,z)\,,\ (x,y)\in\partial\Omega_0, 0\le z-b(x,y)\le h(t,x,y))\}$, inflow/outflow or periodic boundary conditions (for example).

Note that the question of existence and uniqueness of solutions to the Boundary Value Problem (BVP) above is difficult and precisely answered only in a few specific situations, for instance see~\cite{abergel-bona-1992,allain-1985b,beale-nishida-1985,shibata-shimizu-2007} for Newtonian viscous fluids and~\cite{lemeur-2011} for non-Newtonian fluids. 
In this work, we simply {\it assume} that specifying the BCs as above (plus initial conditions) allows one to precisely determine one solution (at least) to the bulk equations.
On the other hand, the fact that we restrict to reduced models that only capture long-wave oscillations of the free surface is also an implicit assumption about the regularity of the solutions to our models (reduced or not).

In the next sections, it will be possible to derive simplified equations that are verified by (univoque) approximations of the solutions to the BVP above in the limit $\e\to0$
only over fixed time ranges $T$, of course.
Then, for the sake of clarity, we rewrite the system of equations using adimensional variables 
that are functions of the non-dimensional {\it scaled} coordinates $(\tilde t, \tilde x,\tilde y,\tilde z)=(t/T,x/L,y/L,z/L)$, which we next abusively still write $(t,x,y,z)$. 
On noting that $p$ and $\Tb$ both have the dimension $(L/T)^2$, we normalize the gravity constant as $gT^2/L$ and obtain adimensional bulk equations which we abusively still write
\begin{equation}
\label{eq:navier-stokes1-adim}
\left\lbrace
\begin{aligned}
D_t \bu & = -\grad p +  \sum_{i=x,y,z} (\partial_x T_{ix}+\partial_y T_{iy}+ \partial_z T_{iz}) \be_i - \f \,, 
\\
\div\bu & = 0 \,,
\end{aligned}
\right.
\end{equation}
without tilde, and that are complemented with the boundary conditions
\begin{subequations}
\begin{alignat}{2}
\label{eq:nopenetration2}
(\bu_H\cdot\grad_H)b 
& =  u_z 
\,,\ && \mbox{for }z=b(x,y),
\\
\label{eq:friction2}
\Tb\bn-\left((\Tb\bn)\cdot\bn\right)\bn - k \left(\bu-(\bu\cdot\bn)\bn\right)  
& = 0\,,\
 && \mbox{for }z=b(x,y),
\\
\label{eq:kinematic2}
\partial_t h + (\bu_H\cdot\grad_H)(b+h)
& = u_z \,,\
 && \mbox{for }z=b(x,y)+ h(t,x,y),
\\
\label{eq:tension2} 
-p\bN+\Tb\bN + \gamma\div(\bN)\bN 
& = 0\,,\
 && \mbox{for }z=b(x,y)+ h(t,x,y).
\end{alignat}
\end{subequations}
where $k$ and $\gamma$ have been scaled by $L/T$ and $L^3/T^2$ respectively.

We obtain simplified equations by successive approximations of the solutions to~(\ref{eq:navier-stokes1-adim}--\ref{eq:nopenetration2}--\ref{eq:friction2}--\ref{eq:kinematic2}--\ref{eq:tension2}) following~\cite{gerbeau-perthame-2001,marche-2007}, first for generic Navier-Stokes equations in Section~\ref{sec:thin}, then specifically for many rheologies in various limit regimes in Sections~\ref{sec:newtonian} (Newtonian fluids),~\ref{sec:power} (power-law fluids) and~\ref{sec:viscoelastic} (Oldroyd-B fluids).
Solutions to those simplified equations shall indeed (formally) approximate solutions to~(\ref{eq:navier-stokes1-adim}--\ref{eq:nopenetration2}--\ref{eq:friction2}--\ref{eq:kinematic2}--\ref{eq:tension2}) insofar as they satisfy the full system of equations plus {\it small} error terms.

\section{A generic long-wave thin-layer framework for free-surface gravity flows}
\label{sec:thin}

A usual assumption to formally reduce the thin-layer flow equations for long-wave oscillations of the free-surface is small topography variations.
Extensions to the case of an arbitrary topography are nontrivial, see e.g.~\cite{bouchut-westdickenberg-2004}.
We thus next content ourself with the quite general case of a topography slowly varying around a plane inclined by a constant angle $\Theta$.
More precisely, we choose in~\eqref{eq:navier-stokes1-adim}
$$
\f = (f_x\equiv+g\sin\Theta,f_y\equiv0,f_z\equiv-g\cos\Theta)
$$
(Navier-Stokes equations are Galilean frame-invariant and rotation is a Galilean change of frame) and we assume 
\beq
\label{H1}
{\rm (H1)}:\ \grad_H b=O(\e) \,.
\eeq

It is also natural to assume bounded horizontal velocities $\bu_H=O(1)$, at $z=b+h$ in particular, and $\div_H\bu_H=O(1)$.
Successive implications can next be ``naturally'' derived following a procedure similar to that in~\cite{gerbeau-perthame-2001,marche-2007}, up to the introduction of an additional assumption (i.e. $k\bu_H|_{z=b}=O(\e)$).

\begin{enumerate}

 \item 

We consider first the mass continuity equation $\div\bu = 0$ with BC~\eqref{eq:nopenetration2}
$$ u_z = \bu_H|_{z=b}\cdot\grad_H b - \int_b^z \div_H\bu_H \,. $$

Using $\bu_H=O(1)$ at $z=b$, this yields $u_z = O(\e)$, and $\grad_H h=O(\e)$ by~\eqref{eq:kinematic2}. 
In fact, assuming that only long-wave oscillations of the free surface matter at first order in $\e$, i.e. $\grad_H h=O(\e)$, is {\it equivalent} to $u_z = O(\e)$ as long as ${\rm (H1)}$ holds, by the formula $u_z = \bu_H|_{z=b+h}\cdot\grad_H (b+h)+\pds{t}h + \int_{b+h}^z \div_H\bu_H$ combining~\eqref{eq:kinematic2} and $\div\bu = 0$.\newline
Reciprocally, the mass continuity equation 
and the BCs~(\ref{eq:nopenetration2},\ref{eq:kinematic2}) are satisfied up to error terms of order $O(\e^{a,a+1,a+1})$, respectively (with $a>0$), if $h,u_z,\bu_H$ are replaced with approximations up to error terms of order $O(\e^{a+1,a+1,a})$, respectively (recall $h\sim\e$), provided of course $\div_H\bu_H$ and $\partial_t h$ are also approximated up to $O(\e^a)$ and $O(\e^{a+1})$ (like $\bu_H$ and $h$).

 \item 

Second, using $u_z = O(\e)$, one infers from the momentum equation projected along $\be_z$ that
\beq
\label{eq:vertical_momentum}
\pds{z}p = f_z +  (\pds{z}T_{zz} + \div_H\Tb_{Hz}) + O(\e) 
\eeq
must be satisfied by $p,T_{zz},\Tb_{Hz}$, as well as by approximations of $p,T_{zz},\Tb_{Hz}$ up to errors $O(\e^{2,2,1})$, respectively (assuming $\pds{z}p,\pds{z}T_{zz},\div_H\Tb_{Hz}$ are approximated up to $O(\e^{1,1,1})$).

Moreover, using $\grad_H h=O(\e)$ and $\grad_H b=O(\e)$, 
one infers $\div\bN(t,x,y)= -\Delta_H (b+h) + O(\e^3)$, so the BC~\eqref{eq:tension2} 
rewrites (recall $\gamma\sim1$ is a constant)
\beq
\label{eq0}
p|_{z=b+h} = -\gamma \Delta_H (b+h) +  (T_{zz} - \Tb_{Hz}\cdot\grad_H(b+h)) + O(\e^3) \,
\eeq
which is satisfied as well by approximations of $h,p,T_{zz},\Tb_{Hz}$ up to errors of order $O(\e^{3,3,3,2})$ (at $z=b+h$ at least for $p,T_{zz},\Tb_{Hz}$, provided $\Delta_H h, \grad_H h$ are approximated as well as $h$).

Since~\eqref{eq0} still holds with $O(\e^3)$ replaced by $O(\e^2)$, one can use~\eqref{eq0} consistently with~\eqref{eq:vertical_momentum} in order to obtain (after integration)
\beq
\label{pression1}
p = f_z(z-(b+h)) - \gamma \Delta_H (b+h) +  T_{zz} - \div_H\int_z^{b+h}\Tb_{Hz} + O(\e^2) \,. 
\eeq
Like~\eqref{eq:vertical_momentum},~\eqref{eq0} with $O(\e^3)$ replaced by $O(\e^2)$ and~\eqref{pression1} are satisfied by approximations of $p,T_{zz},\Tb_{Hz}$ up to errors of order $O(\e^{2,2,1})$, and an approximation of $h$ up to 
$O(\e^2)$,
which in turn requires approximations of $u_z,\bu_H$ up to errors of order $O(\e^{2,1})$ (recall the first item).


 \item 

Third, it is natural to assume $\Delta_H (b+h) = \div_H(\grad_H b + \grad_H h) = O(\e)$ on the one hand, and 
$\max(\Tb_{HH},T_{zz},\Tb_{Hz})|_{z=b+h} = O(1)$ on the other hand. One then gets 
\beq
\label{bc1}
\Tb_{Hz}|_{z=b+h} = (\Tb_{HH}-T_{zz}\I)\grad_H(b+h) + O(\e^2) 
\eeq
from BC~\eqref{eq:tension2} 
with~\eqref{eq0} (even with $O(\e^3)$ replaced by $O(\e^2)$ in~\eqref{eq0}) 
plus the scaling
\beq
\label{eq1}
\Tb_{Hz}|_{z=b+h} = O(\e) \,. 
\eeq


Moreover, from~\eqref{eq:friction2} and $\grad_H b=O(\e)$, we get at $z=b$ using~\eqref{eq:tangent}
\begin{multline*}
\Tb_{Hz}\cdot(\grad_H b)^\bot - (\grad_H b)^\bot\cdot\Tb_{HH}\grad_H b 
\\
= k \bu_H\cdot(\grad_H b)^\bot (1+O(\e^2))\,,
\\
(1-O(\e^2))\Tb_{Hz}\cdot\grad_H b - \grad_H b\cdot\Tb_{HH}\grad_H b + |\grad_H b|^2 T_{zz} 
\\
= k (\bu_H\cdot\grad_H b+u_z|\grad_H b|^2) (1+O(\e^2))\,,
\end{multline*}
thus one obtains, with~\eqref{eq:nopenetration2} and, if $\grad_H b\neq0$, 
$$\Tb_{Hz}= \frac{\Tb_{Hz}\cdot\grad_H b}{|\grad_H b|^2}\grad_H b +\frac{\Tb_{Hz}\cdot(\grad_H b)^\bot}{|\grad_H b|^2}(\grad_H b)^\bot \,,$$
\beq
\label{bc2}
 \Tb_{Hz}|_{z=b} =  (\Tb_{HH}-T_{zz}\I) \grad_H b + k\bu_H(1 + O(\e^2)) + |\Tb_{Hz}| O(\e^2) \,.
\eeq
If $\grad_H b=0$, $\Tb_{Hz}|_{z=b} = k\bu_H$ straightforwardly follows from~\eqref{eq:friction2}, so~\eqref{bc2} remains true.

To proceed, although assuming $\max(\Tb_{HH},T_{zz},\Tb_{Hz})|_{z=b} = O(1)$ seems natural, it is clear that one shall also need an assumption for the scaling of $k\bu_H|_{z=b}$, which is less natural.

Note already that the approximations~(\ref{bc1}--\ref{bc2}) of the BCs~(\ref{eq:tension2},\ref{eq:friction2}) are also satisfied by approximations of $\Tb_{HH}$, $T_{zz}$, $\Tb_{Hz}$ up to error terms of order $O(\e^{1,1,2})$ at least at $z=b,b+h$, as long as~\eqref{eq0} holds (possibly with $O(\e^3)$ replaced by $O(\e^2)$) and $k\bu_H|_{z=b}$ is approximated up to an error $O(\e^2)$.
At this stage, this is in contrast with~\eqref{pression1} which requires approximations of $T_{zz},\Tb_{Hz}$ up to $O(\e^{2,1})$. 
But recall that $\Tb$ is still to be connected to the other unknown variables through additional equations 
once we will have fixed the rheology of the fluid. In fact, most of the work in the sequel
will consist in deriving approximations of the stresses that are coherent with approximations of the other variables. (Anticipating the case of Newtonian fluids, we shall for instance be able to use approximations of $\Tb_{HH},T_{zz},\Tb_{Hz}$ up to $O(\e^{2,2,2})$.)

 \item 

Last, whatever the rheology, we will have recourse to the additional (but usual) assumption
\beq
{\rm (H2)}: k\bu_H|_{z=b}=O(\e)\,. 
\eeq
Indeed, we then next get $ \Tb_{Hz}|_{z=b}=O(\e)$ from~\eqref{bc2}, so $\Tb_{Hz}=O(\e)$ from~\eqref{eq1} and finally 
\beq
\label{eq3}
\pds{z}\Tb_{Hz} \equiv D_t\bu_H + \grad_Hp - \div_H\Tb_{HH} - \f_H = O(1) \,,
\eeq
a ``horizontal'' projection~\eqref{eq3} of the momentum equation that makes sense as a cornerstone for model reduction (all terms are bounded).
Note that~\eqref{eq3} is also satisfied up to additional error terms $O(\e)$ by approximations of $\Tb_{Hz}$, $\bu_H$, $p$, $\Tb_{HH}$ up to $O(\e^{2,1,1,1})$, assuming as usually that $\pds{z}\Tb_{Hz}$, $D_t\bu_H$, $\grad_Hp$, $\div_H\Tb_{HH}$ are then all approximated up to $O(\e)$.


\end{enumerate}

So far, in Section~\ref{sec:thin}, we have (formally) established relations that are necessarily satisfied by smooth solutions to the BVP~(\ref{eq:navier-stokes1-adim}--\ref{eq:nopenetration2}--\ref{eq:friction2}--\ref{eq:kinematic2}--\ref{eq:tension2}) as $h\sim\e\to0$ under the scaling assumptions $\rm (H1)$ and $\rm (H2)$.
Moreover, using ``natural'' assumptions (such that all the terms $\bu_H$, $u_z$, $p$, $\Tb_{HH}$, $\Tb_{Hz}$, $T_{zz}$ that appear in the equations of the initial BVP remain bounded, in particular)
we have also obtained that $u_z,p-T_{zz},\Tb_{Hz}$ are small terms of order $O(\e)$ everywhere in the domain.
The observations above are thus natural candidates for the construction of a {\it reduced model} whose solution coincide with an approximation
\begin{multline}
\label{approximation}
 (h^0,\bu_H^0,u_z^0,p^0-T_{zz}^0,\Tb_{HH}^0-T_{zz}^0\I,\Tb_{Hz}^0)\\
 = (h,\bu_H,u_z,p-T_{zz},\Tb_{HH}-T_{zz}\I,\Tb_{Hz}) + O(\e^{2,1,2,2,1,2}) 
\end{multline}
of a solution to the initial BVP when $\e\to0$.
Such a reduced model could read e.g.
\beq
\label{eq:reduced1}
\left\lbrace
\begin{aligned}
D_t\bu_H^0 + \grad_H (p^0-T_{zz}^0) - \div_H (\Tb_{HH}^0-T_{zz}^0\I) - \pds{z}\Tb_{Hz}^0 - \f_H = 0 \,,
\\
\pds{z}(p^0-T_{zz}^0) - \div_H\Tb_{Hz}^0 - f_z = 0 \,,
\\
\div_H \bu_H^0 + \pds{z}u_z^0 = 0 \,,
\\
(\bu_H^0\cdot\grad_H)b- u_z^0 |_{z=b} = 0 \,,
\\
k\bu_H^0 + (\Tb_{HH}^0-T_{zz}^0\I) \grad_H b - \Tb_{Hz}^0|_{z=b} = 0 \,,
\\
\partial_t h^0 + (\bu_H^0\cdot\grad_H)(b+h^0) - u_z^0|_{z=b+h^0} = 0 \,,
\\
\gamma \Delta_H (b+h^0) + \Tb_{Hz}^0\cdot\grad_H(b+h^0) + (p^0-T_{zz}^0)|_{z=b+h^0}= 0 \,,
\\
(\Tb_{HH}^0-T_{zz}^0\I)\grad_H(b+h^0) - \Tb_{Hz}^0|_{z=b+h^0} = 0 \,,
\\
\end{aligned}
\right.
\eeq
where, in comparison with (\ref{eq:navier-stokes1-adim}--\ref{eq:nopenetration2}--\ref{eq:friction2}--\ref{eq:kinematic2}--\ref{eq:tension2}), the higher-order terms in $\e$ have been forgotten.
Though, a reduced model is {\it coherent} only when a corrected approximation exists 
\begin{multline}
\label{correction}
 (h^0,\bu_H^0,u_z^0,p^0-T_{zz}^0,\Tb_{HH}^0-T_{zz}^0\I,\Tb_{Hz}^0) + O(\e^{2,1,2,2,1,2}) \\ = (h,\bu_H,u_z,p-T_{zz},\Tb_{HH}-T_{zz},\Tb_{Hz})
\end{multline}
that is solution to the initial BVP 
under the same assumptions as those that have been used to define the reduced model
(that is, a solution to~\eqref{eq:reduced1} adequately corrected to solve the initial BVP should yield back~\eqref{eq:vertical_momentum},~\eqref{eq0},~\eqref{bc1},~\eqref{bc2} and \eqref{eq3} with error terms of exactly the same scaling as the one observed above for the solution to the initial BVP). 
Coherence implies that the approximation relationship~\eqref{approximation} may indeed hold (formally as $\e\to0$ for smooth enough solutions) insofar as the correction terms in~\eqref{correction} have the right scaling in order to balance the error terms appearing in the initial BVP due to the definition of the reduced model like~\eqref{eq:reduced1} by truncation of the initial BVP.
Furthermore, similarly to the BVP~(\ref{eq:navier-stokes1-adim}--\ref{eq:nopenetration2}--\ref{eq:friction2}--\ref{eq:kinematic2}--\ref{eq:tension2}),
the system~\eqref{eq:reduced1} is not well-posed without complementing it by equations that come from the rheology of the fluid material that it describes.
Then one expects these equations to be also coherent simplifications of those equations complementing~(\ref{eq:navier-stokes1-adim}--\ref{eq:nopenetration2}--\ref{eq:friction2}--\ref{eq:kinematic2}--\ref{eq:tension2}). 

Before proceeding further to construct reduced models that are coherent for specific rheologies, let us also introduce a generic (widely used) additional manipulation of the initial equations at this stage.
On noting $\bu_H = \frac1h \int_b^{b+h} \bu_H + O(\e)$ whenever e.g. $\partial_z\bu_H=O(1)$ holds (this will often be the case), one in fact often thinks of an approximation $\bu_H^0=\bu_H+O(\e)$ as $\bu_H^0 = \frac1h \int_b^{b+h} \bu_H + O(\e)$.
Then, on noting that acceleration classically rewrites with~(\ref{eq:nopenetration2}--\ref{eq:kinematic2}) using Leibniz rule as
$$
\int_b^{b+h} D_t\bu_H = \pds{t} \int_b^{b+h} \bu_H + \div_H \int_b^{b+h} (\bu_H\otimes\bu_H) \,,
$$
one often considers the ``horizontal'' momentum equations~\eqref{eq3} integrated along $z\in(b,b+h)$ 
\begin{multline}
\label{eq3int}
(\Tb_{Hz}-(\Tb_{HH}-T_{zz}\I)\grad_H(b+h))|_{z=b+h} 
\\
- (\Tb_{Hz}-(\Tb_{HH}-T_{zz}\I)\grad_H b)|_{z=b} + \int_b^{b+h} \f_H
\\
= \pds{t} \int_b^{b+h} \bu_H + \div_H \int_b^{b+h} (\bu_H\otimes\bu_H) + \int_b^{b+h} \grad_H (p-T_{zz}) 
 - \div_H \int_b^{b+h} (\Tb_{HH}-T_{zz}\I)
\end{multline}
where, recalling~\eqref{pression1}, $p-T_{zz}$ can 
be approximated up to an error $O(\e^2)$ using $h$, i.e.
\beq
\label{pression2}
p^0-T_{zz}^0 = f_z(z-(b+h^0)) - \gamma \Delta_H (b+h^0) \,. 
\eeq
Clearly, on using~\eqref{bc2},~\eqref{eq1}, and the integrated continuity equation
\beq
\label{eq:continuity}
\pds{t}h +\div_H \int_b^{b+h} \bu_H = 0
\eeq
to define an approximation $h^0 = h + O(\e^2)$, this is a priori more useful than 
\beq
\label{eq4}
\pds{z}\Tb_{Hz}^0
= D_t\bu_H^0 - \f_H - f_z \grad_H(b+h^0) - \gamma \grad_H \Delta_H (b+h^0) - \div_H(\Tb_{HH}^0-T_{zz}^0\I) \,,
\eeq
as a consequence coherent with~\eqref{eq3}, provided one can close 
\beq
\label{eq5}
\pds{t} \int_b^{b+h} \bu_H + \div_H \int_b^{b+h} (\bu_H\otimes\bu_H) =
- k \bu_H|_{z=b} +  \div_H \int_b^{b+h} (\Tb_{HH}-T_{zz}\I) + h \f_H + O(\e^2) \,, 
\eeq
or the next-order evolution equation for approximate momentum $h^0\bu_H^0 = \int_b^{b+h} \bu_H + O(\e^2)$ 
\begin{multline}
\label{eq5one}
\pds{t} \int_b^{b+h} \bu_H + \div_H \int_b^{b+h} (\bu_H\otimes\bu_H) - \div_H \int_b^{b+h} (\Tb_{HH}-T_{zz}\I)
\\  = - k \bu_H|_{z=b} + h \f_H + h f_z \grad_H(b+h) + h \gamma \grad_H \Delta_H (b+h) + O(\e^3) \,,
\end{multline}
with a coherent approximation of $\int_b^{b+h} (\Tb_{HH}-T_{zz}\I)$ (which typically requires one to vertically integrate the rheological equations), since indeed, whenever $\partial_z\bu_H=O(1)$ holds, it also holds $\frac1h \int_b^{b+h} (\bu_H\otimes\bu_H) = \frac1h \left(\int_b^{b+h}\bu_H\right)\otimes\left(\int_b^{b+h}\bu_H\right) + O(\e)$.
One may also note the latter ``depth-averaged'' approach invokes lower-dimensional variables that depend only on the horizontal coordinates (not $z$), which justifies the label ``reduced model'' (lower-dimensional variables are particularly useful for analytical computations as well as fast numerical simulations).
And it a priori does not seem to necessarily require an explicit approximation of the shear component of the stress $\Tb_{Hz}$ (to close the reduced model at least).
Though, to show the coherence with~\eqref{bc2}, one will in fact still need an expression for $\Tb_{Hz}^0|_{z=b}$ at least; and by the way, this is also often very useful to physical interpretations of the reduced model.



\begin{remark}
Assuming ${\rm (H1)}:\grad_H b = O(\e)$ proved directly connected to our goal of modelling only long-waves, and shows up for instance through the dimension reduction in the reduced model~\eqref{eq:reduced1} (where the variable $T_{zz}$ is not autonomous anymore).
On the contrary, assuming ${\rm (H2)}:k\bu_H|_{z=b}=O(\e)$ is less intuitive 
although it is straightforwardly connected with the useful scaling $\Tb_{Hz}=O(\e)$ (a consequence of $k\bu_H|_{z=b}=O(\e)$ through $\Tb_{Hz}|_{z=b}=k\bu_H|_{z=b}+O(\e)$). 
The latter assumption is key, in fact, to get formal simplifications like~\eqref{eq5} or~\eqref{eq5one}.
Yet, this hypothesis may of course not be true in a number of flows~!
Then one should either use the {\it full} system of equations or another reduced model than the ones derived in the present work (then derived with another strategy). 
\end{remark}

%

\section{Application to Newtonian fluids}
\label{sec:newtonian}

  Internal stresses in Newtonian fluids are defined, after rescaling, with a Reynolds number $\Re$ 
  \beq\label{newtonianstress}
  \Tb = \begin{pmatrix}
	\Tb_{HH} & \Tb_{Hz} \\ \Tb_{Hz}^T & T_{zz}
	\end{pmatrix}
  = \frac1\Re \begin{pmatrix}
	2\mathbf{D}_H(\bu_H) & \pds{z}\bu_H+\grad_H u_z \\ (\pds{z}\bu_H+\grad_H u_z)^T & 2\pds{z}u_z
	\end{pmatrix} \,.
  \eeq
  Without further assumption than $\rm (H1)-(H2)$, one simply obtains (from $u_z=O(\e)$)
  $$
  \Tb = \frac1\Re \begin{pmatrix}
	2\mathbf{D}_H(\bu_H) & \pds{z}\bu_H+O(\e) \\ (\pds{z}\bu_H+O(\e))^T & -2\div_H\bu_H
	\end{pmatrix} \,.
  $$
 Then, to derive a closed reduced model invoking coherent approximations of the stresses (such that $\Tb_{Hz}=O(\e)$ in particular), one needs further assumptions.
 Depending on the treatment of $k\bu_H|_{z=b}=O(\e)$, one can in fact obtain different reduced models in the limit $\e\to0$.

 \subsection{The inertial regime} 

  If we specify ${\rm (H2)}$ as 
  \beq
   {\rm (H2a)}: k\sim\e
  \eeq 
  and, for the scaling of $\Tb$ in~\eqref{newtonianstress} to be compatible with 
  the relations 
  in Section~\ref{sec:thin}, further assume 
  \beq
   {\rm (H3)}: \Re\sim\e^{-1}\,,  \qquad  \text{ and } \qquad {\rm (H4)}: \pds{z}\bu_H=O(1) \,,
  \eeq
  one first obtains $\Tb_{HH},T_{zz}=O(\e)$ and then the improved scalings $\Tb_{HH}-T_{zz}\I,p-T_{zz}=O(\e)$ using ${\rm (H3)}$.
  Moreover, because of $\rm (H4)$, a non-degenerate approximation $\bu^0_H = \bu_H + O(\e)$ that does not go to zero almost everywhere when $\e\to0$ must have a flat profile ($\partial_z \bu_H^0 = 0$), that is
  \beq\label{eq:stratification}
   \bu_H(t,x,y,z) = \bu^0_H(t,x,y) + O(\e) \,,
  \eeq
  also termed a ``motion by slices''. 
  Last, $\rm (H4)$ 
  suffices to justify the depth-averageing procedure introduced at the end of Section~\ref{sec:thin} for the construction of a reduced model,
  so an approximation $(h^0,\bu^0_H)\approx (h,\bu_H)+O(\e^{2,1})$ may be simply determined as a solution to~(\ref{eq:continuity}--\ref{eq5}) 
  where the higher-order terms $O(\e^2)$ have been neglected, i.e. the system 
  \begin{align}
  \label{newtonian1a}
    \pds{t}h^0 + \div_H(h^0\bu_H^0) & = 0 \,, 
    \\
  \label{newtonian1b}
    \pds{t}(h^0\bu_H^0) + \div_H(h^0\bu_H^0\otimes\bu_H^0) + k\bu_H^0 - h^0\f_H & = 0  \,. 
  \end{align}
  But whereas the 
  solutions to~(\ref{newtonian1a}--\ref{newtonian1b})
  straightforwardly allow one to construct a ``first-order'' approximation 
  $(h^0,\bu_H^0,u_z^0,p^0)=(h,\bu_H,u_z,p)+O(\e^{2,1,2,2})$ 
  that is coherent with the continuity equation, with~\eqref{eq:kinematic2} and with the BCs~(\ref{eq:nopenetration2}--\ref{eq0}) as
  \beq
   \label{eq:reconst}
   u_z^0 = \bu_H^0\cdot\grad_H b + (b-z)\div_H\bu_H^0\,, \quad p^0 = f_z(z-(b+h^0))-\gamma\Delta_H(b+h^0)-\frac2\Re\div_H(\bu_H^0) 
  \eeq
  when assumptions ${\rm (H1-H2a-H3-H4)}$ hold, at this stage, 
  one still cannot compute approximations $\Tb_{Hz}^0=\Tb_{Hz}+O(\e^2)$ 
  and it is thus not clear yet that a solution $(h^0,\bu_H^0)$ to~(\ref{newtonian1a}--\ref{newtonian1b}) also defines a coherent approximation of equations~(\ref{eq:navier-stokes1-adim}--\ref{eq:friction2}--\ref{eq:tension2}) used in the derivation of~(\ref{newtonian1a}--\ref{newtonian1b}).
  (Note that it suffices to show that $\Tb_{Hz}^0|_{z=b+h}=O(\e^2)$ and $\Tb_{Hz}^0|_{z=b}=k\bu_H^0+O(\e^2)$ hold.)

  Fortunately, after an adequate combination of \eqref{newtonian1a} and \eqref{newtonian1b}, we note that it holds 
  \beq
   \label{eq:interm}
   \pds{t}\bu_H^0 + (\bu_H^0\cdot\grad_H)\bu_H^0 + k\bu_H^0/h^0 = \f_H 
  \eeq
  so that the approximation proposed above would indeed be a first-order approximate solution to the horizontal projection of the momentum equation~\eqref{eq:navier-stokes1-adim} (that is~\eqref{eq4} without the higher order terms $O(\e)$) if one could construct $\Tb_{Hz}^0$ such that $\partial_z \Tb_{Hz}^0 = k\bu_H^0/h^0 + O(\e)$. 
  In fact, it is now classical that one can achieve this construction thanks to a so-called parabolic correction~\cite{gerbeau-perthame-2001,marche-2007}.
  The point is to construct an approximation $\tilde\bu_H = \bu^0_H + \bu^1_H$ to $\bu_H$, with $\bu^0_H$ a solution to~(\ref{newtonian1a}--\ref{newtonian1b}) plus possibly higher-order terms $O(\e^2)$, and $\bu^1_H=O(\e)$ a 
  correction such that one can characterize $\Tb_{Hz}^0=\Tb_{Hz}+O(\e^2)$. 
  Plugging the ansatz $\tilde\bu_H = \bu^0_H + \bu^1_H$ for $\bu_H$ in~\eqref{eq4} 
  yields, 
  on noting $\pds{z}\Tb_{Hz} = \frac1\Re (\pdds{z}{z}\bu^1_H+\grad_H\div_H\bu_H^0) + O(\e) \equiv \frac1\Re \pdds{z}{z}\bu^1_H + O(\e) $,
  \beq
  \nonumber 
   \pds{t}\bu_H^0 + (\bu_H^0\cdot\grad_H)\bu_H^0
    = 
    \frac{ \pdds{z}{z}\bu^1_H }\Re + O(\e)
  \eeq
  so that, recalling~\eqref{eq:interm} for $\bu^0_H$ solution to~(\ref{newtonian1a}--\ref{newtonian1b}) plus $O(\e^2)$ terms, the correction must satisfy
  \beq
  \label{eq8}
  \frac1\Re \pdds{z}{z}\bu_H^1 
  = D_t\bu_H^0 - \f_H + O(\e) 
  = - \frac{k}{h^0} \bu_H^0 + O(\e) \,.
  \eeq
  Since furthermore $\Re\sim\e^{-1}$ and $\pds{z}\bu_H=O(1)$ in~\eqref{eq1} imply $\Tb_{Hz}|_{z=b+h}=O(\e^2)$,~\eqref{eq8} requires 
  \beq
  \label{eq9}
   \frac1\Re\pds{z}\bu_H^1= k\bu_H^0 \frac{b+h-z}{h} + O(\e^2)\,.
  \eeq
   Now, the trick is to require $\frac1{h}\int_b^{b+h}\bu_H^1=O(\e^2)$, 
   so one can build a coherent approximation $\tilde\bu_H = \bu^0_H + \bu^1_H = \bu_H + O(\e)$ of the initial BVP 
   with a parabolic correction to $\bu^0_H$ 
  \beq
   \label{eq9one}
    \bu_H^1= \frac{\Re k}{2h} \left( (b+3h/2-z)(z-b-h/2) + h^2/12 \right) \bu_H^0 
  \eeq
  that is simply an explicit function of $\bu_H^0$.
  On the other hand, $\bu_H^0$ can indeed be computed coherently with the second-order approximation of the depth-averaged equation~\eqref{eq5one}
  \begin{multline*}
   \pds{t}\left(h\bu_H^0+\int_b^{b+h}\bu_H^1\right) + \div_H\left(h\bu_H^0\otimes\bu_H^0+\bu_H^0\otimes\int_b^{b+h}\bu_H^1+\int_b^{b+h}\bu_H^1\otimes\bu_H^0\right)
   \\
   = h \f_H + h f_z \grad_H(b+h) + h \gamma \grad_H \Delta_H (b+h) + \frac2\Re \div_H\left( h(\mathbf{D}_H(\bu_H^0)+ \div_H\bu_H^0\I) \right)
   \\
   + \frac2\Re \div_H\left(\int_b^{b+h}\mathbf{D}_H(\bu_H^1)+\left(\int_b^{b+h}\div_H\bu_H^1\right)\I\right) - k\tilde\bu_H|_{z=b}  + O(\e^3)
  \end{multline*}
  using $\int_b^{b+h}\bu_H^1=O(\e^3)$, and $\tilde\bu_H|_{z=b}=\bu_H^0 (1-h\Re k/3)$. 
  With the integrated continuity equation 
  and neglecting $O(\e^3)$ terms, one obtains a closed system of equations for $(h^0 = h+ O(\e^2),\bu_H^0)$ 
  \begin{align}
  \pds{t}h^0 + \div_H( h^0 \bu_H^0 ) = 0 
  \label{newtonian2a}
  \\
  \pds{t}(h^0\bu_H^0) + \div_H( h^0\bu_H^0\otimes\bu_H^0 ) - \left( h^0 \f_H + f_z h^0 \grad_H(b+h^0) \right)
  \label{newtonian2b}
  \\
  = \gamma h^0 \grad_H \Delta_H (b+h^0) - k\bu_H^0 (1-h^0\Re k/3)
   + \frac2\Re \div_H\left(h^0\left(\mathbf{D}_H(\bu_H^0)+\div_H\bu_H^0\I\right)\right) 
  \nonumber
  \end{align}
  that also defines, when $\rm (H1-H2a-H3-H4)$ hold, 
  a coherent reduced model 
  for a first-order approximation $(\tilde h,\tilde\bu_H,\tilde u_z,\tilde p)$ 
  of the true solution $(h,\bu_H,u_z,p)$ of the initial BVP. 
  Indeed, if we set $\tilde h = h^0$, if $\tilde\bu_H$ is defined as above using the parabolic correction~\eqref{eq9one} to $\bu_H^0$, if we define $\tilde u_z=u_z^0$ and $\tilde p=p^0$ like in~\eqref{eq:reconst}, then 
  the ``horizontal'' momentum equation is satisfied up to a first-order error term $O(\e)$ (recall~\eqref{eq8}),
   the continuity equation is satisfied up to an $O(\e)$ error term, 
   the ``vertical'' momentum equation is satisfied up to an error term $O(\e)$ (recall~\eqref{eq:vertical_momentum}), 
   and \eqref{eq:nopenetration2}, 
   \eqref{eq:friction2} 
   \eqref{eq:tension2} as well as~\eqref{eq:kinematic2} are satisfied up to an error $O(\e^2)$.

  Note that when~\eqref{eq:friction} is replaced with pure slip ($k=0$), 
  $\bu_H^1$ does not depend on $z$ anymore. 
  Then, from~\eqref{eq8}, the stronger motion-by-slice $\pds{z}\bu_H = O(\e)$ holds, and coherent simplifications do not need explicit approximations for $\Tb_{Hz}$ 
  provided $\Tb_{Hz}=O(\e^2)$. 
  In particular, an approximation $\tilde\bu_H$ can be straightforwardly defined from the solution $\bu_H^0\equiv\tilde\bu_H$ to~(\ref{newtonian2a}--\ref{newtonian2b}) where $k=0$.
  
  Note also that the case of no-tension boundary conditions is obtained by letting $\gamma\to0$ in the reduced model as well as in the initial BVP. (The limit $\gamma\to0$ commutes with $\e\to0$.)

  The scaling regime in this section is termed inertial because the leading terms are purely inertial and the viscous terms appear only in the correction, as opposed to the viscous regime below. 

%

%
%

 \begin{remark}[Inviscid limit and perfect fluids with shallow water equations]
  
  Note that our scaling implies that as $\e\to0$ the full model (Navier-Stokes equations) formally reduces to the incompressible Euler equations, while the reduced model (the so-called viscous 
  shallow water equations) reduces to the inviscid 
  shallow water equations without friction.
  But if we had considered perfect fluids from the beginning, thus $\Tb=\bzero$, the choice of a motion by slice $\rm (H4)$ is not so much a ``natural assumption'' dictated by the internal stresses.
  This shows not only that various limit procedures do not necessarily commute, but also the importance of choosing adequate 
  dissipation terms at the finest level of modelling (even when these terms are small).
  Otherwise, one encounters such infamous difficulties as the modelling of Reynolds stresses that occur in turbulence modelling.

 \end{remark}

\begin{remark}[Dam break, long waves and vorticity with shallow water equations]
 
    It may be a bit surprising that we derive the shalow water 
    system of equations from Navier-Stokes equations under the assumption of small deformation of the free surface.
    Indeed, shallow water 
    equations have been used numerically with success for a long time to simulate dam breaks, a case that does not seem to agree well with $\grad_H h=O(\e)$.
    But this is consistent with the fact that inviscid shallow water equations can also be obtained as a natural limit for potential ideal flows, see e.g.~\cite{bonneton-lannes-2009}, in regimes where surface waves with a short wavelength compared with the water depth are neglected.
    Now, dam breaks indeed correspond to the case of surface waves with an amplitude of order $h$ similar to the water depth that is small compared with wavelengths of order $L$ similar to a supposedly infinitely-long channel, a particular case of tidal waves~\cite{lamb-1975} where viscosity and vorticity are also neglected.
    Besides, note that the scalings used above to obtain the shallow water 
    equations imply in turn that the vorticity has the scaling
    $
      \curl\bu = (O(\e),O(\e),\omega=\grad_H\wedge\bu_H)
    $
    where $\omega=O(\e)$ must also hold since the vorticity equation $D_t(\curl\bu) = [(\curl\bu)\cdot\grad]\bu + \frac1\Re\Delta(\curl\bu)$ implies
    $$
      \omega\pds{z}\bu_H = O(\e) \quad 
      D_t\omega = \omega\pds{z}u_z + \frac1\Re\Delta \omega = \omega\pds{z}u_z + O(\e^2) \Leftrightarrow
      \pds{t}\omega + \div_H(\omega\bu_H) = O(\e^2)
    $$
    on using $\partial_z\omega=O(1)$.
    So a negligible vorticity is not only a sufficient condition to obtain shallow water 
    equations from Euler equations in some cases~\cite{bonneton-lannes-2009}, it also seems necessary, at least in the cases where the scalings of the previous Section above hold.


\end{remark}


 \subsection{The viscous regime}
 
 Instead of assuming ${\rm (H2a)}$ to achieve ${\rm (H2)}$, one can also look for a regime where $k\sim1$ holds and
 \beq
  {\rm (H2b)}: \bu_H|_{z=b}=O(\e) \,.
 \eeq
 Then, one should still require ${\rm (H3)}: \Re\sim\e^{-1}$ and ${\rm (H4)}: \pds{z}\bu_H=O(1)$ in order to next use the observations of Section~\ref{sec:thin} necessarily satisfied by smooth enough solutions.
 On using ${\rm (H2b)}$ and ${\rm (H4)}$, note that it holds $\bu_H=O(\e)$, which is of course stronger than ${\rm (H2b)}$,
 thus also $u_z=O(\e^2)$ by~\eqref{eq:nopenetration2} and the continuity equation.
 This is at the basis of the {\it viscous} regime, where viscous terms dominate in the momentum conservation. 
 In particular, the latter rewrites (recall~\eqref{eq4})
 \beq
  \label{mmomentumviscousa}
  \frac1\Re\pdds{z}{z}\bu_H=\f_H + O(\e)
 \eeq
 and after using~\eqref{eq1} (in fact only $\pds{z}\bu_H|_{z=b+h}=O(\e)$), we obtain 
 \beq
  \label{eqviscous1a}
  \frac1\Re\pds{z}\bu_H=\frac1\Re\pds{z}\bu_H|_{z=b+h}+\f_H(z-(b+h))+O(\e^2)=\f_H(z-(b+h))+O(\e^2).
 \eeq
 Note that if $\Theta=O(\e)$, this yields $\pdds{z}{z}\bu_H=O(1)$, thus $\pds{z}\bu_H=O(\e)$, and
 \beq
  \label{mmomentumviscousb}
  \frac1\Re \pdds{z}{z}\bu_H = \f_H + f_z \grad_H(b+h) + \gamma \grad_H \Delta_H (b+h) + O(\e^2) \,,
 \eeq
 so finally 
 \eqref{eqviscous1a} with $\f_H$ replaced by $\f_H + f_z \grad_H(b+h) + \gamma \grad_H \Delta_H (b+h)$ and an error $O(\e^3)$
 (recalling that~\eqref{eq1} actually implies $\pds{z}\bu_H|_{z=b+h}=O(\e^2)$).

 We next consider the boundary condition~\eqref{bc2} more carefully, it reads
 \beq
 \label{eqviscous2}
  \frac1\Re \pds{z} \bu_H|_{z=b} = k\bu_H + O(\e^3)
 \eeq
 and thus yields $k\bu_H|_{z=b}= - \f_Hh + O(\e^2) = O(\e) $ in the general case, so finally
 \beq
   \label{eqviscous3}
   \bu_H = \f_H \left( \Re \left( (z-(b+h))^2/2-h^2/2 \right) -h/k \right) + O(\e^2)
 \eeq
 (or $k\bu_H|_{z=b}= - \left(\f_H + f_z \grad_H(b+h) + \gamma \grad_H \Delta_H (b+h) \right)h + O(\e^3)=O(\e^2)$ if $\Theta=O(\e)$, thus 
 the scaling $\bu_H = O(\e^2)$, $u_z = O(\e^3)$, and finally yielding
 \eqref{eqviscous1a} with $\f_H$ replaced by $\f_H + f_z \grad_H(b+h) + \gamma \grad_H \Delta_H (b+h)$ and an error $O(\e^3)$).
  
 Finally, we can derive an autonomous equation for $h$ using
 \beq
  \label{disch}
  \int_b^{b+h} \f_H \left( \Re \left( (z-(b+h))^2/2-h^2/2 \right) -h/k \right) = - \f_H \left( \Re \frac{h^3}3 + \frac{h^2}k \right)
 \eeq
 for an approximation of $\int_b^{b+h}\bu_H$ up to order $O(\e^{3})$ (or $O(\e^{4})$ depending on $\Theta$) in the integrated continuity equation $\pds{t}h+\div_H\int_b^{b+h}\bu_H=0$.
 The solution $h^0$ to 
 \beq
  \label{reynolds}
   \partial_t h^0 - \div_H\left( \f_H \left( \Re \frac{|h^0|^3}3 + \frac{|h^0|^2}k \right) \right) = 0
 \eeq
 (with $\f_H + f_z \grad_H(b+h^0) + \gamma \grad_H \Delta_H (b+h^0)$ instead of $\f_H$ when $\Theta=O(\e)$)
 allows one to define a coherent approximation of the initial BVP as long as $\rm (H1-H2b-H3-H4)$ hold, with $\bu_H^0$ reconstructed from $h^0$ following~\eqref{eqviscous3} 
 (slightly modified when $\Theta=O(\e)$) 
 and $u_z^0,p^0$ reconstructed like in the previous section. The stress terms are also easily reconstructed with $\bu_H^0$.
 
 Note that~\eqref{reynolds} is exactly (2.28) in~\cite{RevModPhys.69.931}, where one also comments on the fact that this reduced model is strongly reminiscent of Reynolds lubrication equation~\cite{bayada-chambat-1986} except that here one has a free-surface condition, so the pressure is known to be hydrostatic, while the boundary $z=h$ is unknown.
 One also compute sometimes higher-order approximations of the discharge~\eqref{disch} as a function of $h$ from the momentum equation, 
 see e.g.~\cite{benney-1966,bresch-noble-2007,boutounet-chupin-noble-vila-2008,fernandez-nieto-noble-vila-2010},
 but the resulting models involve high-order derivatives of $h$ (which is a difficulty for numerical simulations) and the coherence of these approximations is not obvious.

\begin{remark}[About the existence of two limit regimes]
\label{rem:tworegimes}
 
 Like the shallow water 
 equations obtained in the inertial regime, the lubrication equation obtained in the viscous regime also has a number of applications, see e.g.~\cite{RevModPhys.69.931}, but this happens in different situations of course.
 A regime where viscous forces dominate the inertial terms to balance gravity seems to suit better with small-scale slow flows (on short times after the flow initiation and in small domains), 
 when boundary effects are important (and $k$ can be chosen as large as necessary to approximate the no-slip boundary condition obtained in the $k\to\infty$ limit).
 On the contrary, an inertial regime seems to suit better to large-scale fast flows (on long times after the flow initiation and in large -- typically geophysical -- domains),
 when boundary effects can be reduced to a small effective friction condition on a fictitious boundary close to the real boundary inward the fluid (thereby defining a boundary layer with limited amplitude).
 Of course this description is only phenomenological and not quantitatively useful.
 Real flows are the result of particular initial and boundary conditions, and adequately choosing one of the two kinds of reduced models 
 (or none, e.g. when boundary effects are important throughout the domain) seems difficult a priori.

 Taking profit of the possible existence of two established regimes, one might also think of combining them: a viscous thin-layer where boundary effects are well taken into account could be 
 physically meaningful
 as a sublayer beneath an inertial thin-layer. For instance, one may want to construct an interface $z=b+Y$ ($0\le Y\le h$) between the two layers that would define 
 a fictitious free-surface for the former (the continuity equation would still yield an autonomous evolution equation for $Y$ where the viscous regime holds) 
 and a 
 fictitious topography (possibly moving) for the latter. 
 There is nevertheless a difficulty: such a construction would necessarily require the horizontal velocities $\bu_H$ to be discontinuous at the interface (at least in the limit $\e\to0$, which implies that $\pds{z}\bu_H$ hence also $\Tb_{Hz}$ is unbounded close to the interface) and the stresses 
 at the interface to satisfy a friction law of Navier type 
 with a coefficient $k$ to be consistently determined with the size $Y$ of the boundary layer. 
 Now, 
 there seems to exist no easy construction of such a friction law $k$ at a fake interface yet,
 and we leave this difficult problem to future works (whatever the rheology).
 For instance, one strategy may be to find a transition layer with depth $\eta\in(0,h-Y)$, $\eta=o(\e)$, a velocity field $U$ solution to the momentum equations in $z\in(b+Y,b+Y+\eta)$, and a friction coefficient $k\sim\e$ (possibly also a tension $\gamma$) such that when $\e\to0$:
 \begin{itemize}
  \item the limit of $U|_{z=b+Y}$ is a good approximation (at $O(\e^2)$) of the limit of $\bu|_{z=(b+Y)^-}$ which is given by the velocity solution to the viscous regime in $z\in(b,b+Y)$, 
  \item the limit of $U|_{z=b+Y+\eta}$ is a good approximation (at $O(\e)$) of the limit of $\pds{z}\bu,\bu|_{z=(b+Y+\eta)^+}$, which is given by the solution to the inertial regime in $z\in(b+Y,b+h)$ with a friction coefficient $k$,
  \item $\pds{z}U/(k\Re U)|_{z=b+Y+\eta}$ has a limit 
  so that Navier friction law holds at $z=b+Y+\eta$,
  \item normal stresses are continuous at $z=b+Y$ (and one may want to define a tension coefficient $\gamma$ at $z=(b+Y)^+$ to formulate this).
 \end{itemize}
  
\end{remark}
 

\section{Application to purely-viscous non-Newtonian fluids}
\label{sec:power}


  Purely viscous non-Newtonian fluids can be described by a power-law model
  \begin{multline}
  \Tb 
  = \frac{|\mathbf{D}(\bu)|^{n-1}}\Re \begin{pmatrix}
	2\mathbf{D}_H(\bu_H) & \pds{z}\bu_H+\grad_H u_z \\ (\pds{z}\bu_H+\grad_H u_z)^T & 2\pds{z}u_z
	\end{pmatrix}
  \\
  = \frac{|\mathbf{D}(\bu)|^{n-1}}\Re \begin{pmatrix}
	2\mathbf{D}_H(\bu_H) & \pds{z}\bu_H+O(\e) \\ (\pds{z}\bu_H+O(\e))^T & -2\div_H\bu_H
	\end{pmatrix}
  = \begin{pmatrix}
	\Tb_{HH} & \Tb_{Hz} \\ \Tb_{Hz}^T & T_{zz}
	\end{pmatrix}
  \end{multline}
  where internal stresses are nonlinear functions of the strain rate due to the non-constant viscosity
  \begin{multline}
  |\mathbf{D}(\bu)|^{n-1}=(|\mathbf{D}_H(\bu_H)|^2+|\pds{z}\bu_H+\grad_H u_z|^2/2+|\pds{z}u_z|^2)^{(n-1)/2}
    \\ =(|\mathbf{D}_H(\bu_H)|^2+|\pds{z}\bu_H+O(\e)|^2/2+|\div_H\bu_H|^2)^{(n-1)/2}.   
  \end{multline}
  The degenerate constant case $n=1$ has been treated in the previous section.
  The cases $0<n<1$ and $n>1$ are clearly different due to different monotonocity properties of the stresses with respect to the deformation gradient $\Dbu$, see e.g.~\cite{barrett-liu-1994}.
  The limit $n\to0$ is singular: it yields a particular occurence of the Bingham model for viscoplastic fluids with a yield stress $|\Dbu|\neq0\Leftrightarrow|\Tb|>\frac2\Re$.

 \subsection{The inertial regime} 
  
    Let us look for a coherent approximation of the solutions to the initial BVP when $\rm (H1-H2a-H3-H4)$ hold, like in the inertial regime of the Newtonian case,
    so that the observations of Section~\ref{sec:thin} are true (at least formally for smooth enough solutions). 
    In fact, only the internal stresses change in the present purely-viscous non-Newtonian case compared with the Newtonian case, and one can follow the same procedure until the construction of a correction.
   Then, the question is: can we proceed, starting from the nonlinear version of~\eqref{eq9}, viz.
  \beq
  \label{eq9power}
   \frac1\Re (|\mathbf{D}_H(\bu_H^0)|^2+\frac{|\pds{z}\bu_H^1|^2}2+|\div_H\bu_H^0|^2)^{(n-1)/2} \pds{z}\bu_H^1= k\bu_H^0 \frac{b+h-z}{h} + O(\e^{2})\,,
  \eeq
  and define a 
  correction $\bu_H^1=O(\e^2)$ using the same trick as in the Newtonian case, that is $\int_b^{b+h}\bu_H^1=O(\e^3)$ ?
  Notice that this is always possible when $\left|\int_b^{b+h}\pds{z}\bu_H^1\right|=O(\e^2)$. 

\smallskip

  For $n\ge1$ (shear-thickening fluids), let us define the function $\phi_a:x\to(x^2/2+a)^{(n-1)/2}x$ that is one-to-one and onto from $\R_{\ge0}$ to $\R_{\ge0}$, so we can rewrite (in componentwise sense)
  $$ \pds{z}\bu_H^1=\phi^{-1}_a( |\Re k\bu_H^0(b+h-z)/h| ) sg( \Re k\bu_H^0(b+h-z)/h ) $$ 
  as a function of $z$ parametrized by $\bu_H^0$ through $a=|\mathbf{D}_H(\bu_H^0)|^2+|\div_H\bu_H^0|^2$.
  Notice that it holds $0\le \phi^{-1}_a( |\Re k\bu_H^0(b+h-z)/h| )\le |\Re k\bu_H^0(b+h-z)/h|$ for $z\in(,b+h)$,
  and since we could do the computation $\int_b^{b+h}|\Re k\bu_H^0(b+h-z)/h|\:dz=O(\e^2)$ exactly in the Netwonian case, 
  it follows that a 
  correction $\bu_H^1$ such that $\int_b^{b+h}\bu_H^1=O(\e^3)$ can be constructed here.
  One can next construct a first-order approximation $\bu_H^0=\bu_H+O(\e)$ with the 
  solution to~(\ref{newtonian2a}--\ref{newtonian2b}) where
  \begin{itemize}
   \item[(i)] $h^0$ is replaced by\footnote{ 
   This is a function of $\bu_H^0$ and $h^0$ that one may obtain numerically after integration $\int_b^{b+h^0}\cdot\:dz$ by quadrature of the terms inside $\div_H$.
   } 
   $\int_b^{b+h^0}\frac{|\Re k\bu_H^0(b+h^0-z)/h^0|}{\phi^{-1}_a( |\Re k\bu_H^0(b+h^0-z)/h^0|)}dz$ in the viscous terms of the RHS of \eqref{newtonian2b}
   \item[(ii)] the friction term invokes the new value\footnote{ 
   This is only known to be bounded above by $h^0\max_{z\in[b,b+h^0]}\pds{z}\bu_H^1\le\Re k\bu_H^0 h^0$. 
  } of $\bu_H|_{z=b}$ approximated at $O(\e^2)$.
  \end{itemize}
  This straightforwardly defines a coherent approximation 
  insofar as the only equation which is different from the (coherent) Newtonian case is the equilibration of second-order viscous dissipation terms with friction at the bottom boundary, 
  and the 
  correction term above has been constructed on purpose for that coherence to be satisfied.

  Note also that this reduced model obtained for $n\ge1$ seems new to us. 
  Though, it is not very practical (because some terms are implicit) and may not be very useful for applications (because shear-thickening fluids are not very common in nature).

\smallskip

  For $n<1$ (shear-thinning fluids), we are not able to conclude about the 
  correction with the strategy above.
  Instead, let us 
  try to compare the cases $0<n<1$ with the limit case $n\to0$.

  For $n=0$, assuming $|\Dbu|\neq0$ (thus $|\mathbf{D}_H(\bu_H^0)|\neq0$)
  the point is again to 
  solve
  \beq
  \label{eq100}
   \frac{\pds{z}\bu_H^1}{\sqrt{|\mathbf{D}_H(\bu_H^0)|^2+|\pds{z}\bu_H^1|^2/2+|\div_H\bu_H^0|^2}} = \frac{\Re k\bu_H^0}{h}( b+h-z ) + O(\e) \,.
  \eeq
  Now, 
  \eqref{eq100} has real solutions $\pds{z}\bu_H^1$ that are compatible with the scaling implied by $\rm (H1-H4)$
  if, and only if, the condition $0<|\bu_H^0|<\sqrt{2}/(k\Re)$ is satisfied. 

  For $0\le n<1$, one can then construct a 
  correction that satisfies $\int_b^{b+h}\bu_H^1=O(\e^3)$
  and 
  define a coherent approximation of the full model provided $0<|\bu_H^0|<\sqrt{2}/(k\Re)$ (componentwise) and $|\mathbf{D}_H(\bu_H^0)|\neq0$: 
  in that case, the reduced model derived for $n\ge1$ still holds.
  When $n\to1$, the limit of that model still coincides (as expected) with the standard (viscous) shallow water equations, 
  and when $n\to0$, 
  the 
  correction can be computed exactly, 
  it has the profile
  \begin{multline}
  \label{profile0}
    \bu_H^1 = \sqrt{|\mathbf{D}_H(\bu_H^0)|^2+|\div_H\bu_H^0|^2} \frac{h}{\Re k|\bu_H^0|} \Bigg( 2 \sqrt{1-\left(\frac{\Re k|\bu_H^0|}{\sqrt{2}h}( b+h-z )\right)^2} \\
+ \sqrt{1-\left(\frac{\Re k|\bu_H^0|}{\sqrt{2}}\right)^2} + \frac{\sqrt{2}}{\Re k|\bu_H^0|} \mathop{\rm arcsin}\left(\frac{\Re k|\bu_H^0|}{\sqrt{2}}\right) \Bigg)\frac{\bu_H^0}{|\bu_H^0|}
    + O(\e^2)  \,.
  \end{multline}
  Moreover, on noting that $|\mathbf{D}_H(\bu_H^0)|=0\Rightarrow|\mathbf{D}_H(\bu_H)|=0$ (otherwise~\eqref{eq100} leads to a contradiction), the case $|\mathbf{D}_H(\bu_H^0)|=0$ should hold if, and only if, $|\Tb|<\frac2\Re$. 
  So one could also obtain the reduced model starting from a variational inequality instead of~\eqref{eq3int} as full model (see e.g.~\cite{duvaut-lions-1972})
  and get similarly to~\cite{bresch-fernandeznieto-ionescu-vigneaux-2010} with a test function $\bv_H$ 
  \begin{align}
  \label{bingham1a}
  & \pds{t}h^0 + \div_H(h^0\bu_H^0) = 0 \,,
   \\
  \label{bingham1b}
  & \int_\Omega \left( \pds{t}(h^0\bu_H^0) + \div_H(h^0\bu_H^0\otimes\bu_H^0) + k(\bu_H^0+\bu_H^1|_b) \right)\cdot(\bv_H-\bu_H^0)
   \\
  \nonumber
  & \quad + \int_\Omega \frac2\Re \left| \mathbf{D}_H(\bv_H)-\mathbf{D}_H(\bu_H^0) \right|
 \ge 
   \\
  \nonumber
 & \qquad \int_\Omega \left( \frac2\Re \div_H\left( \beta \frac{\mathbf{D}_H(\bu_H^0)+\div_H\bu_H^0\I}{\sqrt{|\mathbf{D}_H(\bu_H^0)|^2+|\div_H\bu_H^0|^2}}\right) \right)\cdot(\bv_H-\bu_H^0)
  \\
  \nonumber
  &\qquad + \int_\Omega \left( \gamma h^0 \grad_H \Delta_H (b+h^0) + h^0 \f_H + f_z h^0 \grad_H(b+h^0) \right)\cdot(\bv_H-\bu_H^0)
  \end{align}
  where, using the explicit profile~\eqref{profile0}, one can compute the friction term
  and the viscosity modification $\beta$. 
  But remember that the latter reduced model (in the limit $n\to0$) breaks down when $|\bu_H^0|>\sqrt{2}/(k\Re)$ while, at the same time, it has no meaning when $|\mathbf{D}_H(\bu_H^0)|=0\Leftrightarrow|\Tb|<\frac2\Re$, which seems contradictory.

\begin{remark}[About viscoplastic non-Newtonian fluids]

The physical pertinence and the modelling of viscoplastic non-Newtonian fluids with a yield stress is still much debated. 
In any case, Bingham law is a cornerstone of the viscoplastic modelling since it allows to mathematically investigate the concept of yield-stress and it is worth discussing. That is why we would also like to mention that the most usual form of Bingham law is not as above, but includes an additional viscous dissipative term, and is often thought as a particular case of the more general Herschel-Bulkley law
$$ 
\Tb = \left(\frac2\Re |\Dbu|^{m}+\Bi\right)\frac{\Dbu}{|\Dbu|} \text{ if } \Dbu \neq \bzero, \text{ then } |\Tb|\ge\Bi  \text{, or } \Dbu = \bzero \Leftrightarrow |\Tb|<\Bi \,,
$$
where this time we have denoted $\Bi$ a yield-stress independent of $\frac2\Re$, the usual adimensional constant for the ratio between the viscous dissipation and inertia. The standard Bingham law coincides with the case $m=1$, while we investigated the case $m\to-\infty$ above when $n=0$. 

For any $m$, the conclusion above needs to be modified as follows, provided one assumes $\Bi\sim\e$ in order to perform our thin-layer reduction procedure. 
As above, one cannot go further than derive a reduced-model for the subdomains of the two-dimensional domain $\Omega$ where $|\mathbf{D}_H(\bu_H^0)|\neq0$ holds. 
And the problem still consists in computing a 
correction 
from a profile solution to
\beq
\label{eq110}
 \left(\frac2\Re |\Dbu|^{m}+\Bi\right) \frac{\pds{z}\bu_H^1}{\sqrt{|\mathbf{D}_H(\bu_H^0)|^2+|\pds{z}\bu_H^1|^2/2+|\div_H\bu_H^0|^2}} = \frac{k\bu_H^0}{h}( b+h-z ) + O(\e^2) \,,
\eeq
a polynomial equation in $|\pds{z}\bu_H^1|$ which unfortunately does not seem to be soluble for any $m>0$. 

Another paradigm in viscoplastic modelling has attracted much attention recently, see e.g.~\cite{jop-forterre-pouliquen-2006}, and we would like to mention it too: a Drucker-Prager yield criterion can replace Von Mises one
\beq
\label{eq111}
\Tb = \left(\frac2\Re |\Dbu|^{m}+p\Bi\right)\frac{\Dbu}{|\Dbu|} \text{ if } \Dbu \neq \bzero, \text{ then } |\Tb|\ge p\Bi  \text{, or } \Dbu = \bzero \Leftrightarrow |\Tb|<p\Bi \,.
\eeq
Note that it is not necessary to assume $\Bi=O(\e)$ then since one already has $p=O(\e)$. 
In particular, when
the viscous component vanishes, the 
correction to the velocity profile should satisfy
\beq
\label{eq112}
 p\Bi 
 \frac{\pds{z}\bu_H^1}{\sqrt{|\mathbf{D}_H(\bu_H^0)|^2+|\pds{z}\bu_H^1|^2/2+|\div_H\bu_H^0|^2}} = \frac{k\bu_H^0}{h}( b+h-z ) + O(\e^2) \,,
\eeq
where we recall~\eqref{pression2} $
p = f_z(z-(b+h)) - \gamma \Delta_H (b+h) + T_{zz} + O(\e^2)$. On noting~\eqref{eq111}, it holds
$$
 p \left( 1 + \Bi \frac{\div_H(\bu_H^0)}{\sqrt{|\mathbf{D}_H(\bu_H^0)|^2+|\pds{z}\bu_H^1|^2/2+|\div_H\bu_H^0|^2}} \right)
 = f_z(z-(b+h)) - \gamma \Delta_H (b+h)  + O(\e^2)
$$
which, plugged into~\eqref{eq112}, yields an algebraic equation for $\pds{z}\bu_H^1$ at any $z\in(b,b+h)$
\beq
\label{eq113} 
 \frac{\Bi\pds{z}\bu_H^1\left(f_z(z-(b+h)) - \gamma \Delta_H (b+h)\right)}{ 
   \Bi\div_H(\bu_H^0) + \sqrt{|\mathbf{D}_H(\bu_H^0)|^2+|\pds{z}\bu_H^1|^2/2+|\div_H\bu_H^0|^2} } = \frac{k\bu_H^0}{h}( b+h-z )  + O(\e^2)\,.
\eeq
In the case $\gamma=0$ (no surface tension), the formula becomes much easier
\begin{multline}\label{eq116} 
 \left(\frac12-\left( \frac{h\Bi f_z}{k\bu_H^0} \right)^2\right) |\pds{z}\bu_H^1|^2
  - 2\Bi\div_H(\bu_H^0)\left( \frac{ h\Bi f_z}{k\bu_H^0} \right)|\pds{z}\bu_H^1|
 \\ + |\mathbf{D}_H(\bu_H^0)|^2+(1-\Bi^2)|\div_H\bu_H^0|^2 + O(\e^2) = 0 
\end{multline}
and one can then also solve explicitly the problem for the 
correction.
So the solution to~\eqref{eq116} allows one to define an admissible velocity correction, and thus also a coherent approximation of the full model through the reduced model, 
as soon as 
the sole requirement
$|\mathbf{D}_H(\bu_H^0)|\neq0\Leftrightarrow|\mathbf{D}_H(\bu_H)|\neq0$ 
is satisfied here (a condition that unfortunately remains difficult to predict or analyze here ; in particular, we are not aware of a simpler reformulation of this model as a variational inequality).

\end{remark}

 \subsection{The viscous regime}
 
  Assuming $\rm (H1-H2b-H3-H4)$ we again follow, for purely viscous non-Newtonian fluids, the same procedure as in the Newtonian case.
  First we obtain a nonlinear version of~\eqref{eqviscous1a}
  \beq
  \label{eqviscous1anonlinear}
   \frac1\Re  ( |\pds{z}\bu_H|^2/2 )^{(n-1)/2} \pds{z} \bu_H 
   \\ = \f_H \left( z-(b+h) \right) + O(\e^2)
  \eeq
  on noting $\mathbf{D}_H(\bu_H)=O(\e)$ (with additional terms to $\f_H$ if $\Theta=O(\e)$).
  With the friction boundary condition at $z=b$, this next yields 
  \beq
  \label{eqviscous2nonlinear}
   \bu_H = \left( \Re 2^{\frac{n-1}2} \bolda \right)^{\frac1n} \left( \frac{(z-(b+h))^{\frac{n+1}n} - h^{\frac{n+1}n}}{\frac{n+1}n} \right) - \bolda \frac{h}k + O(\e^{1+\frac{2/3}n}) \,,
  \eeq
  where $\bolda = \f_H$, or $\bolda =\f_H + f_z\grad_H(b+h) + \gamma \grad_H \Delta_H (b+h)$ if $\Theta=O(\e)$,
  and an autonomous equation for $h^0=h+O(\e^2)$ from the continuity equation $\pds{t}h+\div_H\int_b^{b+h}\bu_H=0$ and the approximation
  \beq
  \label{eqviscous3nonlinear}
    \int_b^{b+h}\bu_H  =  \left( \Re 2^{\frac{n-1}2} \bolda \right)^{\frac1n} \left(\frac{2n+1}{n+1} h^{\frac{2n+1}n}\right)- \bolda \frac{h^2}k + O(\e^{2+\frac{2/3}n})\,.
  \eeq
  This coincides with the viscous limit recently derived in~\cite{fernandez-nieto-noble-vila-2010}, though with another mathetically-inclined viewpoint and a slightly different scaling (the term $h^2/k$ is absent in particular, somehow a no-slip limit $k\to\infty$).
  It holds for all power-law fluids (though, note that the quality of approximation increases with $n$ in the shear-thinning case but decreases in the shear-thickening case).


\section{Application to viscoelastic non-Newtonian fluids}
\label{sec:viscoelastic}

There are numerous models for viscoelastic non-Newtonian fluids, with various definitions of the extra-stress $\str$ in $\Tb=2\eta_s\Dbu+\str$.
We 
concentrate here on one prototypical model among 
{\it differential} constitutive equations for $\str$,
 the Upper-Convected Maxwell (UCM) equations~\cite{bird-curtiss-armstrong-hassager-1987a},
\begin{equation}
\label{eq:ucm}
D_t\str = (\gbu)\str+\str(\gbu)^T + \frac{1}{\lambda}\left(2\eta_p\Dbu-\str\right) \qquad\mbox{in } \D(t),
\end{equation}
where $\lambda$ is interpreted as a characteristic relaxation time for elastic dilute molecules 
and $\eta_p$ as a viscosity. 
There are many extensions to the UCM equations, which one also often writes using the total (kinematic) viscosity $\eta=\eta_s+\eta_p$ and the retardation time $\lambda(1-\theta)\le\lambda$ where $\theta=\eta_p/\eta\in(0,1)$
\begin{equation}
\label{eq:navier-stokes2}
\left\lbrace
\begin{aligned}
D_t\bu & = -\grad p + \div(2\eta(1-\theta)\Dbu) + \div\str + \f & \text{ in } \D(t) \,,
\\
\div\bu & = 0 & \text{ in } \D(t) \,,
\\
\lambda( D_t\str - (\gbu)\str-\str(\gbu)^T ) & = 2\eta\theta\Dbu-\str & \text{ in } \D(t).
\end{aligned}
\right.
\end{equation}
A simple one for instance combines the power-law and the UCM models $\Tb=2\eta_s|\Dbu|^{n-1}\Dbu+\str$, see~\cite{mohamed-reddy-2010}. 
One can also use nonlinear versions of the relaxation term in the right-hand side of~\eqref{eq:ucm}, 
see~\cite{renardy-2000}. 
But~\eqref{eq:ucm} 
already contains the kinematic essence of most {\it differential} constitutive equations (material frame indifference for tensors)
and we postpone the discussion of other models to Remark~\ref{FENEP} (and possible future works).



  To adimensionalize~\eqref{eq:navier-stokes2}, let us 
  introduce the Deborah number $\De = \lambda/T$, and 
$$
  \Tb = \begin{pmatrix}
	\Tb_{HH} & \Tb_{Hz} \\ \Tb_{Hz}^T & T_{zz}
	\end{pmatrix}
  = 
  \frac{1-\theta}\Re \begin{pmatrix}
	2\mathbf{D}_H(\bu_H) & \pds{z}\bu_H+\grad_Hu_z \\ (\pds{z}\bu_H+\grad_Hu_z)^T & 2\pds{z}u_z 
	\end{pmatrix}
  + \begin{pmatrix}
	\str_{HH} & \str_{Hz} \\ \str_{Hz}^T & \tau_{zz}
	\end{pmatrix} 
$$
  so the extra-stress $\str$ satisfies  the non-dimensional UCM equations
  \begin{equation}
  \De \left( D_t\str - (\gbu)\str - \str(\gbu)^T \right) = \frac{2\theta}{\Re}\Dbu-\str \,.
  \label{eq:ucm2}
  \end{equation}

  Note by the way that the cases $\De=O(\e)$ are 
  not the most physically interesting
  because they lead us back to a purely-viscous Newtonian extra-stress at first order of approximation in $\e\to0$.

  In the following, we use the reformulation of~\eqref{eq:ucm2} with the also well-known {\it conformation tensor} variable $\strs = \I + \frac{\De\Re}\theta\str$ 
  solution to an evolution equation using the single scalar parameter $\De$
  \begin{equation}
  \De \left( D_t\strs - (\gbu)\strs - \strs(\gbu)^T \right) = \I-\strs \,.
  \label{eq:ucm3}
  \end{equation}
  The infamous Weissenberg number $\Wi = \De\Re/\theta$ then 
  appears in Navier-Stokes~\eqref{eq:navier-stokes1-adim} through
  \beq
\label{stressvisco}
  \Tb  = 
  \frac{1-\theta}\Re \begin{pmatrix}
	2\mathbf{D}_H(\bu_H) & \pds{z}\bu_H+O(\e) \\ (\pds{z}\bu_H+O(\e))^T & -2\div_H\bu_H
	\end{pmatrix}
  + \frac\theta{\Re\De} 
   \begin{pmatrix}
	\strs_{HH}-\I_H & \strs_{Hz} \\ \strs_{Hz}^T & \sigma_{zz}-1
    \end{pmatrix} 
  \,,
  \eeq
  where, recalling Section~\ref{sec:thin}, we have also used the continuity equation, $h\sim\e$ and $\rm (H1)$ in 
  $$
  \gbu = \begin{pmatrix}
	\grad_H \bu_H & \pds{z} \bu_H \\ (\grad_H u_z)^T & \pds{z}u_z
	\end{pmatrix}
   = \begin{pmatrix}
	\grad_H \bu_H & \pds{z} \bu_H \\ O(\e) & -\div_H\bu_H
	\end{pmatrix} \,.
  $$
  We recall that for physical reasons\footnote{
    More precisely, for the model to be consistent with the usual thermodynamics principles, see e.g.~\cite{wapperom-hulsen-1998-a}.
  }, the conformation tensor should always be \textit{positive-definite}, and indeed remains so as long as it is initially and solutions to~\eqref{eq:ucm3} are smooth enough (see e.g.~\cite{boyaval-lelievre-mangoubi-2009}).

  From now on, recalling Section~\ref{sec:thin}, it is natural to assume that $\strs_{HH}$ and $\sigma_{zz}$ are not only bounded but also have the same scaling.
  Then, on noting that~\eqref{eq:ucm3} reads
  \begin{subequations}
  \beq
  \hspace{-2mm}
  \De \left( D_t\strs_{HH} - (\grad_H\bu_H)\strs_{HH} - \strs_{HH}(\grad_H\bu_H)^T - \strs_{Hz}\otimes\pds{z}\bu_H - \pds{z}\bu_H\otimes\strs_{Hz} \right) = \strs_{HH}-\I 
  \label{eq:ucm3HH}
   \eeq
  \beq
  \De \left( D_t\strs_{Hz} - (\grad_H\bu_H)\strs_{Hz} - \strs_{HH}(\grad_Hu_z) - \strs_{Hz}\pds{z}u_z - \pds{z}\bu_H\sigma_{zz} \right)  = \strs_{Hz} 
  \label{eq:ucm3Hz}
  \eeq
  \beq
  \De \left( D_t\sigma_{zz} - 2\strs_{Hz}\cdot\grad_Hu_z - 2\sigma_{zz}\pds{z}u_z \right)  = \sigma_{zz}-1 \,,
  \label{eq:ucm3zz}
  \eeq
  \end{subequations}
  it stems from~\eqref{eq:ucm3HH} and~\eqref{eq:ucm3Hz}
  that 
  ${\rm (H4)}: \pds{z}\bu_H=O(1)$ is also as natural (for boundedness) in viscoelastic non-Newtonian fluids as in purely viscous (Newtonian and non-Newtonian) fluids. 
  Under $\rm (H4)$, one then obtains with $\De\sim1$ 
  \begin{subequations}
  \beq
  \hspace{-13mm}
  \De \left( D_t\strs_{HH} - (\grad_H\bu_H^0)\strs_{HH} - \strs_{HH}(\grad_H\bu_H^0)^T - \strs_{Hz}\otimes\pds{z}\bu_H^1 - \pds{z}\bu_H^1\otimes\strs_{Hz} \right) 
   = \strs_{HH}-\I + O(\e)
  \label{eq:ucm3HH1}
  \eeq
  \beq
  \De \left( D_t\strs_{Hz} - (\grad_H\bu_H^0)\strs_{Hz} + \strs_{Hz}\div_H\bu_H^0 - \pds{z}\bu_H^1\sigma_{zz} \right)  = \strs_{Hz} + O(\e)
  \label{eq:ucm3Hz1}
  \eeq
  \beq
  \De \left( D_t\sigma_{zz} + 2\sigma_{zz}\div_H\bu_H^0 \right)  = \sigma_{zz}-1  + O(\e)
  \label{eq:ucm3zz1}
  \eeq
  \end{subequations}
  from~(\ref{eq:ucm3HH}--\ref{eq:ucm3Hz}--\ref{eq:ucm3zz}), for any first-order approximation $\bu_H^0=\bu_H+O(\e)$ with a flat profile, 
  possibly corrected 
  by some $\bu_H^1=O(\e)$.
 
  \subsection{The inertial regime}

  Like in the previous cases, we obtain an inertial limit when one specifies $\rm (H2)$ as ${\rm (H2a)}: k\sim\e$.
  On the contrary, to coherently use $ \Tb_{Hz} = O(\e) $ for BCs~\eqref{eq1} and~\eqref{bc2} in Section~\ref{sec:thin} with
  \beq
   \label{shearviscoelastic}
   \Tb_{Hz} = \frac1\Re \left( (1-\theta)\pds{z}\bu_H^1 + \theta \frac1\De\strs_{Hz} \right) \,,
  \eeq
  we should now further assume, in addition to $\rm (H4)$,
  \begin{itemize}
   \item[(i)] either ${\rm (H3)}: \Re\sim\e^{-1}$ like in the Newtonian case, 
   \item[(ii)] or ${\rm (H5a)}: 1-\theta\sim\e$, plus either ${\rm (H6a)}: \strs_{Hz}=O(\e)$ or ${\rm (H6c)}: \De\sim\e^{-1}$, 
   \item[(iii)] or ${\rm (H5b)}: \pds{z}\bu_H=O(\e)$ (which is stronger than $\rm (H4)$) plus either $\rm (H6a)$, or ${\rm (H6b)}: \theta\sim\e$, or ${\rm (H6c)}$. 
  \end{itemize}
  Note that in absence of other 
  assumptions, $\De\sim1$ and $\theta\sim1$ shall be simply taken as constants.
  
\subsubsection{Small internal stresses}

  Under assumptions $\rm (H1-H2a-H4-H3)$, like in the Newtonian case, first-order approximations $(h^0,\bu_H^0)$ solution to~(\ref{newtonian1a}--\ref{newtonian1b}) are not necessarily coherent with BCs~\eqref{eq1} and~\eqref{bc2} and the point is how to approximately compute $\Tb_{Hz}$. 
  Introducing a 
  correction $\bu_H^1$ satisfying
  \beq 
  \label{eq10}
    \frac1\Re \left( (1-\theta)\pds{z}\bu_H^1 + \theta \frac1\De\strs_{Hz} \right) = k\bu_H^0\frac{b+h-z}h + O(\e^2) \,,
  \eeq
  one would then like to coherently replace~(\ref{newtonian1a}--\ref{newtonian1b}) by a reduced model invoking the depth-averaged horizontal momentum equation truncated at order $O(\e^3)$ (so the impact of the 
  correction $\bu_H^1$ on $\bu_H^0$ is coherently taken into account), just like in the Newtonian case, plus simplified UCM equations to close the system. 
  Now, the dissipative terms involving the viscoelastic stress tensor $\Tb$ in the momentum equation can be computed
  using approximations of the UCM system of equations 
  without any explicit reference to the 
  correction $\bu_H^1$ after rewriting~\eqref{eq10}
  \beq
  \label{eq10bis}
   \pds{z}\bu_H^1 = \frac1{1-\theta}\left(\Re k\bu_H^0\frac{b+h-z}h-\theta \frac1\De\strs_{Hz}\right) + O(\e) \,.
  \eeq
  Then a coherent reduced model is obtained as usual after closing the second-order truncation of the horizontal momentum equation, typically using the same trick $\int_b^{b+h}\bu^1_H=O(\e^3)$ as in the Newtonian case. Assuming, for the sake of simplicity,
  $$
   {\rm (H7a)}: \pds{z}\strs_{HH},\pds{z}\sigma_{zz}=O(1) \qquad    {\rm (H7b)}: \pds{z}\strs_{Hz}=O(1)
  $$
  a profile can be computed explicitly from~\eqref{eq10bis} such that $\int_b^{b+h}\bu^1_H=O(\e^3)$ holds 
  \beq
  \label{eq10ter}
   \bu_H^1 = \frac1{1-\theta}\left(\frac{\Re k}{2h}\bu_H^0 \left( \frac{h^2}3-(b+h-z)^2 \right) - \frac{\theta}{2} \frac1\De\strs_{Hz}^0 \left(z-(h+2b)\right) \right) \,,
  \eeq
  and a reduced model coherent at first-order 
  with $\rm (H1-H2a-H4-H3-H7)$ reads
  \begin{subequations}
  \beq
  \pds{t}h^0 + \div_H( h^0 \bu_H^0 ) = 0 
  \label{nonnewtonian1a}
  \eeq
  \begin{multline}
  \label{nonnewtonian1b}
   \pds{t}(h^0\bu_H^0) + \div_H( h^0\bu_H^0\otimes\bu_H^0 )   + k\bu_H^0 \left(1-\frac{\Re}{(1-\theta)}\frac{kh^0}3\right)
  + k \strs_{Hz}^0 \frac{\Re}{(1-\theta)}\frac{\theta(b+h^0)}{2\De}
\\
  = \left( h^0 \f_H + f_z h^0 \grad_H(b+h^0) \right) + \gamma h^0 \grad_H \Delta_H (b+h^0) 
  \\
   + \frac{2(1-\theta)}\Re \div_H\left(h^0\left(\mathbf{D}_H(\bu_H^0)+\div_H\bu_H^0\I\right)\right) 
   + \frac\theta{\Re\De} \div_H \left(h^0(\strs_{HH}^0-\sigma_{zz}^0\I)\right) 
  \end{multline}
  \begin{multline}
  \label{nonnewtonian1c}
  \De \left( \pds{t}(h^0\strs_{HH}^0)+\div_H(h^0\bu_H^0\otimes\strs_{HH}^0) \right) 
  \\ = 
   h^0 \frac{\De}{1-\theta} \left( \strs_{Hz}^0\otimes\left(\frac{\Re k}2\bu_H^0-\theta\frac1\De\strs_{Hz}^0\right)
         - \left(\frac{\Re k}2\bu_H^0-\theta\frac1\De\strs_{Hz}^0\right)\otimes\strs_{Hz}^0 \right)
   \\
    + h^0 \De \left( (\grad_H\bu_H^0)\strs_{HH}^0 + \strs_{HH}^0(\grad_H\bu_H^0)^T \right) + h^0(\strs_{HH}^0-\I)
  \end{multline}
  \begin{multline}
  \label{nonnewtonian1d}
  \De \left( \pds{t}(h^0\strs_{Hz}^0)+\div_H(h^0\bu_H^0\otimes\strs_{Hz}^0) \right)  
\\
   =
  h^0 \frac{\De}{1-\theta} \left(\frac{\Re k}2\bu_H^0-\theta\frac1\De\strs_{Hz}^0\right) \sigma_{zz}^0
  + h^0 \De \left(\grad_H\bu_H^0\strs_{Hz}^0 - \div_H\bu_H^0\strs_{Hz}^0 \right)
  + h^0\strs_{Hz}^0
  \end{multline}
  \beq 
  \label{nonnewtonian1e}
  \De \left( \pds{t}(h^0\sigma_{zz}^0)+\div_H(h^0\bu_H^0\sigma_{zz}^0) \right) = h^0 \De ( 2\sigma_{zz}^0\div_H\bu_H^0 ) + h^0(\sigma_{zz}^0-1)
  \,.
  \eeq
  \end{subequations}

  Note that one retrieves the standard viscous 
  shallow water (our reduced model for the standard Navier-Stokes equations) in the limit $\theta\to0$  (prior or subsequent to $\e\to0$ ; i.e. the two formal limits commute here), plus UCM equations that then become simply enslaved transport equations for a material tensor (without feedback in the momentum equation).
  Otherwise, this seems to be a new model. 
  In particular, it was not identified in our previous work~\cite{bouchut-boyaval-2013} that focused on the case $\theta=1$ (where tangential boundary conditions like friction are a priori useless constraint for the initial BVP) because then, it is not possible to derive an expression for $\pds{z}\bu_H$ (with a link between the shear strain and the shear stress like~\eqref{eq10} one cannot compute a coherent approximation of~\eqref{eq:ucm3Hz1} like~\eqref{nonnewtonian1d}). 
  Unfortunately, the limit $\theta\to1$ after $\e\to0$ is unclear, but we next assume $\theta=1+O(\e)$ without even assuming $\Re\sim\e^{-1}$ then, 
  and will next be able to derive a limit model provided $\strs_{Hz}/\De=O(\e)$.

  
\subsubsection{Small viscous internal stresses}

  Under assumptions $\rm (H1-H2a-H5a-H6a)$, a non-vanishing first-order approximation of $\strs_{Hz}$ 
  can be coherently constructed from~\eqref{eq:ucm3Hz} only if ${\rm (H5b)}$ holds (since $\sigma_{zz} = O(\e)$ is impossible, having as equilibrium value $1$ by~\eqref{eq:ucm3zz} as long as $\De$ remains bounded), which is on the other hand not coherent with~\eqref{eq10} 
  and the horizontal momentum equation unless $\bu_H^0=0$. 
  That is why we only consider $\rm (H1-H2a-H5a-H6c)$, 
  plus ${\rm (H7)}$ for the sake of simplicity, which leads to the reduced model
  \begin{subequations}
  \beq
  \pds{t}h^0 + \div_H( h^0 \bu_H^0 ) = 0 
  \label{nonnewtonian3a}
  \eeq
  \begin{multline}
  \label{nonnewtonian3b}
  \pds{t}(h^0\bu_H^0) + \div_H( h^0\bu_H^0\otimes\bu_H^0 )   + k\bu_H^0 \left(1-\frac{\Re}{(1-\theta)}\frac{kh^0}3\right) + k \strs_{Hz}^0 \frac{\Re}{(1-\theta)}\frac{\theta(b+h^0)}{2\De}\\
   = \left( h^0 \f_H + f_z h^0 \grad_H(b+h^0) \right) + \gamma h^0 \grad_H \Delta_H (b+h^0) 
  \\ + \frac{2(1-\theta)}\Re \div_H\left(h^0\left(\mathbf{D}_H(\bu_H^0)+\div_H\bu_H^0\I\right)\right) 
   + \frac\theta{\Re\De} \div_H \left(h^0(\strs_{HH}^0-\sigma_{zz}^0\I)\right) 
  \end{multline}
  \begin{multline}
  \label{nonnewtonian3c}
  \pds{t}(h^0\strs_{HH}^0)+\div_H(h^0\bu_H^0\otimes\strs_{HH}^0) = 
  h^0 \left( (\grad_H\bu_H^0)\strs_{HH}^0 + \strs_{HH}^0(\grad_H\bu_H^0)^T  \right)
  \\
   + h^0 \frac1{1-\theta} \left( \strs_{Hz}^0\otimes\left(\frac{\Re k}2\bu_H^0-\theta\frac1\De\strs_{Hz}^0\right)
         - \left(\frac{\Re k}2\bu_H^0-\theta\frac1\De\strs_{Hz}^0\right)\otimes\strs_{Hz}^0 
   \right)
  \end{multline}
  \begin{multline}
  \label{nonnewtonian3d}
   \pds{t}(h^0\strs_{Hz}^0)+\div_H(h^0\bu_H^0\otimes\strs_{Hz}^0) = h^0(\grad_H\bu_H^0)\strs_{Hz}^0 - h^0\strs_{Hz}^0\div_H\bu_H^0 
    \\ + h^0 \frac1{1-\theta} \left(\frac{\Re k}2\bu_H^0-\theta\frac1\De\strs_{Hz}^0\right) \sigma_{zz}^0
  \end{multline}
  \beq
  \label{nonnewtonian3e}
   \pds{t}(h^0\sigma_{zz}^0)+\div_H(h^0\bu_H^0\sigma_{zz}^0) = 2h^0\sigma_{zz}^0\div_H\bu_H^0 
  \eeq
  \end{subequations}
  whose solutions are coherent with first-order approximations of the initial BVP when $\rm (H1-H2a-H5a-H6c-H7ab)$ indeed holds.

  The latter model {\rm (6.13)} 
  formally coincides with 
  the limit $1/\De\to0$, 
  $\theta\to1$ (provided $k/(1-\theta)$ and $1/\De(1-\theta)$ remain bounded) of {\rm (6.12)}, 
  some kind of ``High-Weissenberg limit'' 
  (where the UCM model suffers from deficiencies, see e.g.~\cite{boyaval-lelievre-mangoubi-2009}, and Remark~\ref{FENEP} for repair suggestions).


\subsubsection{Small viscous internal shear stresses}

  Under assumptions $\rm (H1-H2a-H5b)$, the motion-by-slice is stronger than the usual one, 
  which thus further restricts a priori the regimes of validity of a possible reduced model (even if the reduced model had solutions beyond the regime of validity of our assumptions, such solutions would not necessarily define coherent approximations of the initial BVP).
  It implies 
  \beq \label{veryflat}
   \bu_H(t,x,y,z) = \bu_H^0(t,x,y) + O(\e^2)
  \eeq
  so that the correction $\bu_H^1$ to 
  $\bu_H^0$ is of higher-order than usual ones 
  and does not show up in the horizontal momentum equation if, on the other hand, the extra-stress terms can be computed coherently.
  Now, under $\rm (H1-H2a-H5b-H6a-H7ab)$ -- $\rm (H7)$ for the sake of simplicity -- one indeed obtains the following reduced model coherent with first-order approximations of the initial BVP
  \begin{subequations}
  \beq
  \pds{t}h^0 + \div_H( h^0 \bu_H^0 ) = 0 
  \label{nonnewtonian4a}
  \eeq
  \begin{multline}
  \label{nonnewtonian4b}
  \pds{t}(h^0\bu_H^0) + \div_H( h^0\bu_H^0\otimes\bu_H^0 ) + k\bu_H^0
  \\
   = \left( h^0 \f_H + f_z h^0 \grad_H(b+h^0) \right) + \gamma h^0 \grad_H \Delta_H (b+h^0) 
  \\ + \frac{2(1-\theta)}\Re \div_H\left(h^0\left(\mathbf{D}_H(\bu_H^0)+\div_H\bu_H^0\I\right)\right) 
   + \frac\theta{\Re\De} \div_H \left(h^0(\strs_{HH}^0-\sigma_{zz}^0\I)\right) 
  \end{multline}
  \begin{multline}
  \label{nonnewtonian4c}
  \De \left( \pds{t}(h^0\strs_{HH}^0)+\div_H(h^0\bu_H^0\otimes\strs_{HH}^0) \right)
   \\
   = h^0 \De \left( (\grad_H\bu_H^0)\strs_{HH}^0 + \strs_{HH}^0(\grad_H\bu_H^0)^T \right) + h^0(\strs_{HH}^0-\I)
  \end{multline}
  \begin{multline}
   \label{nonnewtonian4d}
  \De \left(    \pds{t}(h^0\strs_{Hz}^0)+\div_H(h^0\bu_H^0\otimes\strs_{Hz}^0)  \right)
\\
   =   
   h^0 \De \strs_{HH}^0\left(\grad_H(\bu_H^0\cdot\grad_Hb)+(\div_H\bu_H^0)\grad_Hb-\frac12h^0\grad_H\div_H\bu_H^0\right)
\\
  + h^0 \frac{\De}{1-\theta} \left(\frac{\Re k}2\bu_H^0-\theta\frac1\De\strs_{Hz}^0\right)\sigma_{zz}^0
  + h^0 \De \left( (\grad_H\bu_H^0)\strs_{Hz}^0 - h^0 \strs_{Hz}^0\div_H\bu_H^0 \right) + h^0 \strs_{Hz}^0
  \end{multline}
  \beq 
  \label{nonnewtonian4e}
  \De \left( \pds{t}(h^0\sigma_{zz}^0)+\div_H(h^0\bu_H^0\sigma_{zz}^0) \right) = h^0 \De ( 2\sigma_{zz}^0\div_H\bu_H^0 ) + h^0(\sigma_{zz}^0-1)
  \,.
  \eeq
  \end{subequations}
  where, contrary to~(\ref{nonnewtonian1a}--\ref{nonnewtonian1b}--\ref{nonnewtonian1c}--\ref{nonnewtonian1d}--\ref{nonnewtonian1e}) or its ``High-Weissenberg limit''~(\ref{nonnewtonian3a}--\ref{nonnewtonian3b}--\ref{nonnewtonian3c}--\ref{nonnewtonian3d}--\ref{nonnewtonian3e}), the shear component $\strs_{Hz}$ of the viscoelastic stress decouples from the autonomous system of equations~(\ref{nonnewtonian4a}--\ref{nonnewtonian4b}--\ref{nonnewtonian4c}--\ref{nonnewtonian4e}) and is simply computed as a post-processed solution to~\eqref{nonnewtonian4d} enslaved through $\bu_H^0$. (In~\eqref{nonnewtonian4d}, we have used~\eqref{eq10} for the vertical derivative of the horizontal velocity, and the 
  approximate vertical velocity $u_z^0 = u_z+O(\e^2)$ 
  reconstructed from $\bu_H^0$, the continuity equation and the impermeability condition at the bottom excatly like in the Newtonian case, so~\eqref{nonnewtonian4d} is coherent with a first-order approximation $\strs_{Hz}^0=\strs_{Hz}+O(\e^2)$.)

  The latter reduced model {\rm (6.15)} is exactly the viscous two-dimensional extension of the one-dimensional model derived in~\cite{bouchut-boyaval-2013} for the case $\theta=1$, $k=0$.
  The case $k=0$ (pure-slip boundary condition at bottom close to the topography) for $\theta\in[0,1)$ is straightforwardly recovered by taking the limit $k\to0$ in the system above. 
  One cannot compute directly the case $\theta\to1$ ; we refer to~\cite{bouchut-boyaval-2013} for the singular case $\theta=1$ (with $k=0$, where $1^0$ the computation of $\strs_{Hz}^0$ from the horizontal momentum equation supplemented with the bottom boundary conditon, and $2^0$ the computation of an approximation of $\pds{z}\bu_H^1$ from an equivalent to~\eqref{nonnewtonian4d}, are modified).

  When one assumes $\rm (H6b-H7ab)$ in addition to $\rm (H1-H2a-H5b)$, one straightforwardly obtains the same autonomous system of equations as in the reduced model with $\rm (H6a)$, that is~(\ref{nonnewtonian4a}--\ref{nonnewtonian4b}--\ref{nonnewtonian4c}--\ref{nonnewtonian4e}).
  But although it defines a coherent first-order approximation without even assuming any scaling for $\strs$
  (a coefficient $\theta$ of the whole tensor is then only responsible for the small scale), 
  a first-order approximation $\strs_{Hz}^0=\strs_{Hz}+O(\e)$ of a shear component that is not smaller than the other components of the viscoelastic stress tensor would then be different, and i.e. solve 
  \beq
  \De \left( D_t\strs_{Hz}^0 - (\grad_H\bu_H^0)\strs_{Hz}^0 + \strs_{Hz}^0\div_H\bu_H^0 \right)  = \strs_{Hz}^0 \,.
  \label{eq:ucm3Hz1trunc}
  \eeq

  Last, $\rm (H1-H2a-H5b-H6c-H7ab)$ yields the following reduced model coherent with a first-order ``High-Weissenberg-limit'' approximation of the initial BVP
  \begin{subequations}
  \beq
  \pds{t}h^0 + \div_H( h^0 \bu_H^0 ) = 0 
  \label{nonnewtonian5a}
  \eeq
  \begin{multline}
  \label{nonnewtonian5b}
  \pds{t}(h^0\bu_H^0) + \div_H( h^0\bu_H^0\otimes\bu_H^0 )   + k\bu_H^0 
 \\
  = \left( h^0 \f_H + f_z h^0 \grad_H(b+h^0) \right) + \gamma h^0 \grad_H \Delta_H (b+h^0) 
  \\ + \frac{2(1-\theta)}\Re \div_H\left(h^0\left(\mathbf{D}_H(\bu_H^0)+\div_H\bu_H^0\I\right)\right) 
   + \frac\theta{\Re\De} \div_H \left(h^0(\strs_{HH}^0-\sigma_{zz}^0\I)\right) 
  \end{multline}
  \beq
  \label{nonnewtonian5c}
  \pds{t}(h^0\strs_{HH}^0)+\div_H(h^0\bu_H^0\otimes\strs_{HH}^0) = 
  h^0(\grad_H\bu_H^0)\strs_{HH}^0 + \strs_{HH}^0(\grad_H\bu_H^0)^T 
  \eeq
  \beq
  \label{nonnewtonian5d}
   \pds{t}(h^0\strs_{Hz}^0)+\div_H(h^0\bu_H^0\otimes\strs_{Hz}^0) = (\grad_H\bu_H^0)\strs_{Hz}^0 - \strs_{Hz}^0\div_H\bu_H^0 
  \eeq
  \beq
  \label{nonnewtonian5e}
   \pds{t}(h^0\sigma_{zz}^0)+\div_H(h^0\bu_H^0\sigma_{zz}^0) = 2h^0\sigma_{zz}^0\div_H\bu_H^0  \,.
  \eeq
  \end{subequations}

  The various latter models obtained under assumption $\rm (H5b)$ cannot be easily linked to any other one.
  But it is remarkable that in any case, no correction to the flat profile is necessary under assumption $\rm (H5b)$ (even if a profile can be reconstructed afterwards from~\eqref{eq10ter}), whereas the presence of purely (Newtonian) viscous forces is in turn hardly seen but in dissipation terms when one enforces $\rm (H5b)$ instead of $\rm (H5a)$.
  Furthermore, requiring the velocity to have a flat profile~\eqref{veryflat} is thus a priori a very strong limit for the applicability of our reduced models to real flows.
  This may however be particularly interesting for the cases where the normal stress differences are large, since the stress~\eqref{stressvisco} then reads
  \beq
  \label{stressvisco1}
  \Tb  = 
  \frac{1-\theta}\Re \begin{pmatrix}
	2\mathbf{D}_H(\bu_H) & O(\e) \\ O(\e) & -2\div_H\bu_H
	\end{pmatrix}
  + \frac\theta{\Re\De} 
   \begin{pmatrix}
	\strs_{HH}-\I_H & \strs_{Hz} \\ \strs_{Hz}^T & \sigma_{zz}-1
    \end{pmatrix} 
  \,,
  \eeq
  where either ${\rm (H6b)}: \theta\sim\e$, eor  ${\rm (H8)}: \strs_{HH} = \I + O(\e)\,,\ \sigma_{zz}=1+O(\e)\,,\ \strs_{Hz}=O(\e^2)$ (a stronger assumption necessary when starting with ${\rm (H6a)}: \strs_{Hz}\sim\e$) or ${\rm (H6c)}: \De\sim\e^{-1}$ holds but the viscous stretch need not be scaled even though viscoelastic components are always small.  

  To conclude this section, note that even though some reduced models have been identified in the High-Weissenberg limit regime ${\rm (H6c)}: \De\sim\e^{-1}$ where already the model is questionable, we have obtained otherwise two main reduced models -- the autonomous systems of equations~(\ref{nonnewtonian1a}--\ref{nonnewtonian1b}--\ref{nonnewtonian1c}--\ref{nonnewtonian1d}--\ref{nonnewtonian1e}) and~(\ref{nonnewtonian5a}--\ref{nonnewtonian5b}--\ref{nonnewtonian5c}--\ref{nonnewtonian5e}) -- 
  whose solutions define coherent approximations of the initial BVP in physically sensible regimes.
  It could be interesting to numerically simulate the first one, which has not been done yet to our knowledge, and maybe compare it with two-dimenionsal extensions of the solutions to the second model computed in~\cite{bouchut-boyaval-2013}.
  Note in particular that shear effects are then not necessarily small in comparison with elongational/compression effects, 
  which was a problem for the applicability of the second reduced model to real (often sheared~!) flows already noted in~\cite{bouchut-boyaval-2013}.

 \subsection{The viscous regime}

 Assuming $\rm (H1-H2b-H4)$, we proceed for the viscous limit of viscoelastic fluids as usual. 
 We specify $\rm (H2)$ as ${\rm (H2b)}: \bu_H|_{z=b}=O(\e)$ and next require $\Tb_{Hz}=O(\e)$ as above in the inertial case, in addition to ${\rm (H4)}: \pds{z}\bu_H=O(1)$.
 Recall also that the flow is necessary slow here ($\bu_H=O(\e)$) and one obtains from the momentum balance
 \beq
  \label{eqviscousOB1a}
   \frac1\Re \left( (1-\theta)\pds{z}\bu_H + \theta \frac1\De\strs_{Hz} \right) = \f_H(z-(b+h)) 
   + O(\e^2)
 \eeq
 after using 
 $\Tb_{Hz}|_{z=b+h}=O(\e^2)$ and $\int_z^{b+h} \div_H(\Tb_{HH}-T_{zz})=O(\e^2)$.


 Assuming ${\rm (H3)}:\Re\sim\e^{-1}$ plus $\rm (H7)$ for the sake of simplicity in addition to $\rm (H1-H2b-H4)$ 
 (and of course $\De\sim1$, $\theta\sim1$ as long as nothing different is precised for these adimensional numbers)
 leads to a reduced model that is an autonomous system of equations for $(h^0,\strs_{Hz}^0,\sigma_{zz}^0)$
 \begin{subequations}
  \beq 
  \label{eq:ucm3HzviscousOB0}
    \pds{t}h^0 + \frac1{1-\theta} \div_H \left( \frac{\Re}6\f_H |h^0|^3 - \theta \frac1\De \strs_{Hz}^0 \frac{|h^0|^2}2 \right) = 0
  \eeq
   \beq
  \De \left( \partial_t(h^0\strs_{Hz}^0)  \right)  
    = h^0 \strs_{Hz}^0 + h^0 \frac{\De}{1-\theta} \left( \frac{\Re}2\f_H - \theta \frac1\De \strs_{Hz}^0 \right) \sigma_{zz}^0 \,,
  \label{eq:ucm3HzviscousOB1}
  \eeq
  \beq
  \De \left( \partial_t(h^0\sigma_{zz}^0) \right) = h^0 (\sigma_{zz}^0-1) \,, 
  \label{eq:ucm3zzviscousOB1}
  \eeq
  \end{subequations}
  where the discharge in the continuity equation is computed from~\eqref{eqviscousOB1a} and
  \beq
   \label{vel1}
   \bu_H = \frac1{1-\theta} \left( \frac{\Re}2\f_H \left(z-(b+h)\right)^2 - \theta \frac1\De \strs_{Hz}^0 (z-b) \right) + O(\e^2) \,,
  \eeq
  and the longitudinal (horizontal) stress components are obtained by the post-processing
  \begin{multline}
  \De \left( \partial_t(h^0\strs_{HH}^0) \right) = \strs_{HH}^0-\I  
   \\ + h^0 \frac1{1-\theta} \De \left( 
     \strs_{Hz}^0 \otimes \left( \frac{\Re}2\f_H - \theta \frac1\De \strs_{Hz}^0 \right) 
      \left( \frac{\Re}2\f_H - \theta \frac1\De \strs_{Hz}^0 \right) \otimes \strs_{Hz}^0 
   \right) \,.
  \label{eq:ucm3HHviscousOB1}
  \end{multline}
  If furthermore $\Theta=O(\e)$, then the same reduced model hold, but it yields coherent approximations as long as $\pds{z}\bu_H=O(\e)$ and $\strs=O(\e)$ hold, on noting the starting point
  \beq
  \label{eqviscousOB2a}
   \frac1\Re \left( (1-\theta)\pds{z}\bu_H + \theta \frac1\De\strs_{Hz} \right) = (\f_H- f_z \grad_H(b+h) - \gamma \grad_H \Delta_H (b+h))(z-(b+h)) 
  + O(\e^3) \,.
  \eeq

   
  Assuming ${\rm (H5a)}: 1-\theta\sim\e$ and ${\rm (H6a)}: \strs_{Hz}=O(\e)$ again requires $\pds{z}\bu_H=O(\e)$ 
  and cannot be coherent, so we consider ${\rm (H5a)}$ with ${\rm (H6c)}: \De\sim\e^{-1}$ only, which leads to $\sigma_{zz} = 1 + O(\e)$ constant (equal to physical equilibrium) and a reduced model consisting of the limits of~\eqref{eq:ucm3HzviscousOB0} and~\eqref{eq:ucm3HzviscousOB1} as $1/\De\to0$ (with obvious specificities if $\Theta\sim\e$, and~\eqref{eq:ucm3HHviscousOB1} for post-processing $\strs_{HH}^0$ only).

  Last, assuming ${\rm (H5b)}: \pds{z}\bu_H=O(\e)$ and ${\rm (H6b)}: \theta=O(\e)$ leads to the same reduced model as the first one above (and is coherent under the more restricitive regime where $O(\e^2)$ is replaced by $O(\e^3)$ in~\eqref{vel1}), while ${\rm (H5b)}$ and ${\rm (H6c)}$ gives the same as the second one above.

  All these systems seem new to us: 
  other viscous limits of non-Newtonian viscoelastic fluid models have already been derived, but on assuming different scalings, see e.g.~\cite{bayada-chupin-martin-2007,bayada-chupin-grec-2009,bayada-chupin-martin-2010} ($\De\sim\e$). 

  \begin{remark}[Nonlinear differential constitutive equations and HWNP]
  \label{FENEP} 

  The most used variations of the UCM model are {\it nonlinear} modifications of these differential constitutive equations, for instance the FENE-P model where the extra-stress reads
  $\str = \frac\theta{\De\Re} \left( \frac{\strs}{1-\tr\strs/b} - \I \right)$, 
  $b>0$ is a new parameter such that $0\le\tr\strs\le b$ is prerserved by smooth time evolutions of the flow, and the conformation tensor $\strs$ is solution to the nonlinear equation
  \begin{equation}
  \De \left( D_t\strs - (\gbu)\strs - \strs(\gbu)^T \right) = \I-\frac{\strs}{1-\tr\strs/b} \,.
  \label{eq:fenep}
  \end{equation}

  One nice feature of these nonlinear versions is that they usually impose such constraints as $0\le\tr\strs\le b$ which are believed to alleviate the deficiencies of the UCM model (High-Weissenberg-Number Problems or HWNP in short) in the ``High-Weissenberg limit'' (at least, well-posedness has sometimes been shown for smooth flows, see e.g.~\cite{masmoudi-2011}).

  Furthemore, for most of them, reduced models are easily derived from the UCM reduced models above as long as one does not use ${\rm (H6a)}: \strs_{Hz}=O(\e)$.
  It suffices to multiply the last term on the right by $\frac1{1-\tr\strs/b}$, which is indeed never small,
  \begin{itemize}
   \item in~\eqref{nonnewtonian1b},~\eqref{nonnewtonian1c} and~\eqref{nonnewtonian1e} under $\rm (H1-H2a-H4-H3-H7)$,
   \item in~\eqref{nonnewtonian5b},~\eqref{nonnewtonian5c} and~\eqref{nonnewtonian5e} under $\rm (H1-H2a-H4-H5b-H6b-H7)$.
  \end{itemize}
  On the contrary, since $\frac1{1-\tr\strs/b}$ can become arbitrary large when $\tr\strs\to b$, this is not only incompatible with ${\rm (H6a)}: \strs_{Hz}=O(\e)$, but also requires additional assumptions in the case ${\rm (H6c)}: \De\sim\e^{-1}$ (thus not treated here).

  Another way to avoid HWNP is to assume $\De\sim\e$ like in e.g.~\cite{bayada-chupin-martin-2007,bayada-chupin-grec-2009,bayada-chupin-martin-2010}~!
  Then, one cannot expect strong viscoelastic influences on the flow, of course.
  Though, this scaling may be enough for some applications, and we would like to mention that it has recently raised interesting new persepctives:
  a new approach to formal model reduction combining micro and macro scales~\cite{narbona-reina-bresch-2010} that is indeed consistent with a Newtonian behaviour in the limit $\De\to0$.

  \end{remark}

 
%
%
%

 \section{Conclusion}

 We have defined a mathematical framework that allows, for many fluids (i.e. many rheologies), to derive {\it coherent} long-wave thin-layer approximations of free-surface Navier-Stokes flows driven by gravity above smoothly varying topographies.
 Most reduced models derived herein were already known, and the shallow water equations in particular have already proved useful in the numerical simulation of dam-breaks for instance.
 On the other hand, the models for viscoelastic fluids seem to have been much less explored, and some of those derived herein seem new to us.
 Of course, the question how well they model real flows is still to be answered.
 This could be investigated numerically in future works (letting alone their well-posedness and the mathematical control of their distance to the full model) following the same path as in our previous work~\cite{bouchut-boyaval-2013} where we considered the viscoelastic model without friction nor surface tension in a one-dimensional (fast) inertial flow regime.
 Note by the way that the present work also answers important questions concerning the ability of long-wave thin-layer reduced models at describing viscoelastic fluids in {\it sheared} (as opposed to purely extensional) inertial flow regimes that were raised in~\cite{bouchut-boyaval-2013}.
 Last, the non-Newtonian viscous fluids with power-law models and their viscoplastic limit have attracted much attention recently, in particular with a view to modelling avalanches and debris flows. Indeed, such complex flows seem to require complex rheologies, possibly with a yield stress and nonlinear effects,
 while many difficulties have been encountered so far as concerns the modelling of fluid/solid transitions.
 We hope that the unified framework derived herein will help 
 characterize features  essential to long-wave thin-layer flow modelling, and evaluate future models for mud flows and landslides in particular. 

\section*{Acknowledgments}

Thanks to Enrique Fernandez-Nieto for fruitful discussions. 

\providecommand{\bysame}{\leavevmode\hbox to3em{\hrulefill}\thinspace}
\providecommand{\MR}{\relax\ifhmode\unskip\space\fi MR }
\providecommand{\href}[2]{#2}


\begin{thebibliography}{10}

\bibitem{abergel-bona-1992}
F.~Abergel and J.~L. Bona, \emph{A mathematical theory for viscous,
  free-surface flows over a perturbed plane}, Arch. Rational Mech. Anal.
  \textbf{118} (1992), no.~1, 71--93. \MR{1151927 (92k:35213)}

\bibitem{achdou-pironneau-valentin-1998}
Y.~Achdou, O.~Pironneau, and F.~Valentin, \emph{Effective boundary conditions
  for laminar flows over periodic rough boundaries}, J. Comp. Phys.
  \textbf{147} (1998), 187--218.

\bibitem{allain-1985b}
G.~Allain, \emph{Small-time existence for the {N}avier-{S}tokes equations with
  a free surface and surface tension}, Free boundary problems: application and
  theory, {V}ol.\ {IV} ({M}aubuisson, 1984), Res. Notes in Math., vol. 121,
  Pitman, Boston, MA, 1985, pp.~355--364. \MR{903286 (88h:35093)}

\bibitem{ancey-2007}
Christophe Ancey, \emph{Plasticity and geophysical flows: A review}, Journal of
  Non-Newtonian Fluid Mechanics \textbf{142} (2007), no.~1-3, 4 -- 35, In {\it
  Viscoplastic fluids: From theory to application.}

\bibitem{ancey-cochard-2009}
Christophe Ancey and Steve Cochard, \emph{The dam-break problem for
  herschel–bulkley viscoplastic fluids down steep flumes}, Journal of
  Non-Newtonian Fluid Mechanics \textbf{158} (2009), no.~1-3, 18 -- 35,
  Visco-plastic fluids: From theory to application.

\bibitem{barnes-hutton-walters-1989}
H.A. Barnes, J.F. Hutton, and K.F.R.S. Walters, \emph{An introduction to
  rheology}, 1st ed., Elsevier Science Publisher, Amsterdam, The Netherlands,
  1989.

\bibitem{barrett-liu-1994}
John~W. Barrett and W.~B. Liu, \emph{Quasi-norm error bounds for the finite
  element approximation of a non-{N}ewtonian flow}, Numer. Math. \textbf{68}
  (1994), no.~4, 437--456. \MR{1301740 (95h:65078)}

\bibitem{basson-gerardvaret-2008}
Arnaud Basson and David G{\'e}rard-Varet, \emph{Wall laws for fluid flows at a
  boundary with random roughness}, Comm. Pure Appl. Math. \textbf{61} (2008),
  no.~7, 941--987. \MR{2410410 (2009h:76055)}

\bibitem{bayada-chupin-grec-2009}
G.~Bayada, L.~Chupin, and B.~Grec, \emph{{Viscoelastic fluids in thin domains:
  a mathematical proof}}, Asymptotic analysis \textbf{64} (2009), no.~3,
  185--211.

\bibitem{bayada-chupin-martin-2007}
G.~Bayada, L.~Chupin, and S.~Martin, \emph{{Viscoelastic fluids in a thin
  domain}}, Quarterly of Applied Mathematics \textbf{65} (2007), no.~4,
  625--652.

\bibitem{bayada-chupin-martin-2010}
\bysame, \emph{{Viscoelastic Fluids in a Thin Domain: A Mathematical Study for
  a Non-Newtonian Lubrication Problem}}, Mathematical Modeling, Simulation,
  Visualization and e-Learning (D.~Konat\'e, ed.), Springer, 2010,
  pp.~315--321.

\bibitem{bayada-chambat-1986}
Guy Bayada and Mich\`ele Chambat, \emph{The transition between the stokes
  equations and the reynolds equation: A mathematical proof}, Applied
  Mathematics and Optimization \textbf{14} (1986), no.~1, 73--93.
  \MR{MR0826853}

\bibitem{beale-nishida-1985}
J.~Thomas Beale and Takaaki Nishida, \emph{Large-time behavior of viscous
  surface waves}, Recent topics in nonlinear {PDE}, {II} ({S}endai, 1984),
  North-Holland Math. Stud., vol. 128, North-Holland, Amsterdam, 1985, Lect.
  Notes Numer. Appl. Anal. (Vol. 8), pp.~1--14. \MR{882925 (88f:35121)}

\bibitem{benney-1966}
D.~J. Benney, \emph{Long waves on liquid films}, J. Math. Phys. \textbf{45}
  (1966), 50--155.

\bibitem{bird-curtiss-armstrong-hassager-1987a}
R.~B. Bird, C.~F. Curtiss, R.~C. Armstrong, and O.~Hassager, \emph{Dynamics of
  polymeric liquids}, vol. 1: {F}luid {M}echanics, John Wiley \& Sons, New
  York, 1987.

\bibitem{bird-curtiss-armstrong-hassager-1987b}
\bysame, \emph{Dynamics of polymeric liquids}, vol. 2: {K}inetic {T}heory, John
  Wiley \& Sons, New York, 1987.

\bibitem{bonneton-lannes-2009}
P.~Bonneton and D.~Lannes, \emph{Derivation of asymptotic two-dimensional
  time-dependent equations or surface water wave propagation}, Physics of
  Fluids \textbf{2009} (21), 016601.

\bibitem{bouchut-boyaval-2013}
Fran{\c c}ois Bouchut and S{\'e}bastien Boyaval, \emph{{A new model for shallow
  viscoelastic fluids}}, M3AS \textbf{23}, no.~8, 1479--1526 (2013).

\bibitem{bouchut-westdickenberg-2004}
Francois Bouchut and Michael Westdickenberg, \emph{Gravity driven shallow water
  models for arbitrary topography}, Commun. Math. Sci. \textbf{2} (2004),
  no.~3, 359--389. \MR{MR2118849 (2005m:76026)}

\bibitem{boyaval-lelievre-mangoubi-2009}
S{\'e}bastien Boyaval, Tony Leli{\`e}vre, and Claude Mangoubi,
  \emph{Free-energy-dissipative schemes for the {O}ldroyd-{B} model}, M2AN
  Math. Model. Numer. Anal. \textbf{43} (2009), no.~3, 523--561. \MR{2536248
  (2010k:65197)}

\bibitem{bresch-fernandeznieto-ionescu-vigneaux-2010}
D.~Bresch, E.~D. Fern\`andez-Nieto, I.~R. Ionescu, and P.~Vigneaux,
  \emph{Augmented lagrangian method and compressible visco-plastic flows:
  Applications to shallow dense avalanches}, New Directions in Mathematical
  Fluid Mechanics (Giovanni~P. Galdi, John~G. Heywood, Rolf Rannacher,
  Andrei~V. Fursikov, and Vladislav~V. Pukhnachev, eds.), Advances in
  Mathematical Fluid Mechanics, Birkh\"auser Basel, 2010,
  10.1007/978-3-0346-0152-8\_4, pp.~57--89.

\bibitem{bresch-noble-2007}
Didier Bresch and Pascal Noble, \emph{Mathematical justification of a shallow
  water model}, Methods Appl. Anal. \textbf{27} (2007), no.~2, 87--118.

\bibitem{ruyerquil-manneville-2000}
{C. Ruyer-Quil} and {P. Manneville}, \emph{Improved modeling of flows down
  inclined planes}, Eur. Phys. J. B \textbf{15} (2000), no.~2, 357--369.

\bibitem{chupin-2009a}
Laurent Chupin, \emph{The {FENE} viscoelastic model and thin film flows}, C. R.
  Math. Acad. Sci. Paris \textbf{347} (2009), no.~17--18, 1041--1046.

\bibitem{chupin-2009b}
\bysame, \emph{The {FENE} viscoelastic model and thin film flows}, Methods and
  Applications of Analysis \textbf{16} (2009), no.~1, 217--262.

\bibitem{matar-kraster-2009}
R.~V. Craster and O.~K. Matar, \emph{Dynamics and stability of thin liquid
  films}, Rev. Mod. Phys. \textbf{81} (2009), 1131--1198.

\bibitem{degond-lemou-picasso-2002}
P.~Degond, M.~Lemou, and M.~Picasso, \emph{Viscoelastic fluid models derived
  from kinetic equations for polymers}, SIAM J. Appl. Math. \textbf{62} (2002),
  no.~5, 1501--1519 (electronic). \MR{1918565 (2003j:82074)}

\bibitem{duvaut-lions-1972}
G.~Duvaut and J.-L. Lions, \emph{Les in\'equations en m\'ecanique et en
  physique}, Dunod, 1972.

\bibitem{fernandez-nieto-noble-vila-2010}
E.~D. Fern{\'a}ndez-Nieto, P.~Noble, and J.-P. Vila, \emph{Shallow water
  equations for non newtonian fluids}, Journal of Non-Newtonian Fluid Mechanics
  \textbf{165} (2010), no.~13--14, 712--732.

\bibitem{gerardvaret-2009}
David G{\'e}rard-Varet, \emph{The navier wall law at a boundary with random
  roughness}, Communications in Mathematical Physics \textbf{286} (2009),
  81--110 (English).

\bibitem{gerbeau-perthame-2001}
Jean-Fr\'ed\'eric Gerbeau and Beno\^it Perthame, \emph{Derivation of viscous
  {S}aint-{V}enant system for laminar shallow water ; numerical validation},
  Discrete and continuous dynamical system Series B \textbf{1} (2001), no.~1,
  89--102.

\bibitem{jager-mikelic-2003}
Willi J{\"a}ger and Andro Mikeli{\'c}, \emph{Couette flows over a rough
  boundary and drag reduction}, Comm. Math. Phys. \textbf{232} (2003), no.~3,
  429--455. \MR{1952473 (2003j:76025)}

\bibitem{jop-forterre-pouliquen-2006}
Pierre Jop, Yo\"el Forterre, and Olivier Pouliquen, \emph{A constitutive law
  for dense granular flows}, Nature \textbf{441} (2006), no.~7094, 727--730.

\bibitem{lamb-1975}
Lamb, \emph{Hydrodynamics}, Cambridge University Press, 1975.

\bibitem{boutounet-chupin-noble-vila-2008}
P.~Noble J.-P.~Vila M.~Boutounet, L.~Chupin, \emph{Shallow water equations for
  newtonian fluids over arbitrary topographies}, Comm Math Sci \textbf{6}
  (2008), no.~1, 29--55.

\bibitem{marche-2007}
F.~Marche, \emph{{Derivation of a new two-dimensional viscous shallow water
  model with varying topography, bottom friction and capillary effects}},
  European Journal of Mechanics-B/Fluids \textbf{26} (2007), no.~1, 49--63.

\bibitem{masmoudi-2011}
N.~Masmoudi, \emph{Global existence of weak solutions to macroscopic models of polymeric flows},
  Journal de Math\'ematiques Pures et Appliqu\'ees \textbf{96} (2011), no.~5, 502--520.

\bibitem{lemeur-2011}
Herv\'e~Le Meur, \emph{Well-posedness of surface wave equations above a
  viscoelastic fluid}, Journal of Mathematical Fluid Mechanics \textbf{13}
  (2011), no.~4, 481--514.

\bibitem{mohamed-reddy-2010}
H.~B.~H. Mohamed and B.~D. Reddy, \emph{Some properties of models for
  generalized {O}ldroyd-{B} fluids}, Internat. J. Engrg. Sci. \textbf{48}
  (2010), no.~11, 1470--1480, Special Issue in Honor of K.R. Rajagopal.
  \MR{2760996}

\bibitem{narbona-reina-bresch-2010}
G.~Narbona-Reina and D.~Bresch, \emph{On a shallow water model for
  non-newtonian fluids}, Numerical Mathematics and Advanced Applications 2009
  (Gunilla Kreiss, Per L{\"o}tstedt, Axel M\o{a}lqvist, and Maya Neytcheva,
  eds.), Springer Berlin Heidelberg, 2010, 10.1007/978-3-642-11795-4 74,
  pp.~693--701.

\bibitem{RevModPhys.69.931}
Alexander Oron, Stephen~H. Davis, and S.~George Bankoff, \emph{Long-scale
  evolution of thin liquid films}, Rev. Mod. Phys. \textbf{69} (1997), no.~3,
  931--980.

\bibitem{renardy-2000}
M.~Renardy, \emph{Mathematical analysis of viscoelastic flows}, CBMS-NSF
  Conference Series in Applied Mathematics, vol.~73, SIAM, 2000.

\bibitem{shibata-shimizu-2007}
Yoshihiro Shibata and Senjo Shimizu, \emph{On a free boundary problem for the
  {N}avier-{S}tokes equations}, Differential Integral Equations \textbf{20}
  (2007), no.~3, 241--276. \MR{2293985 (2008a:35211)}

\bibitem{wapperom-hulsen-1998-a}
Peter Wapperom and Martien~A. Hulsen, \emph{Thermodynamics of viscoelastic
  fluids: the temperature equation}, J. Rheol. \textbf{42} (1998), no.~5,
  999--1019.

\end{thebibliography}
\end{document}